# ASYMPTOTIC NORMALITY OF PLUG-IN LEVEL SET ESTIMATES

By David M. Mason[1] and Wolfgang Polonik[2]

*University of Delaware and University of California, Davis*

We establish the asymptotic normality of the $G$-measure of the symmetric difference between the level set and a plug-in-type estimator of it formed by replacing the density in the definition of the level set by a kernel density estimator. Our proof will highlight the efficacy of Poissonization methods in the treatment of large sample theory problems of this kind.

**1. Introduction.** Let $f$ be a Lebesgue density on $\mathbb{R}^d$, $d \geq 1$. Define the level set of $f$ at level $c \geq 0$ as

$$C(c) = \{x : f(x) \geq c\}.$$

In this paper we are concerned with the estimation of $C(c)$ for a given level $c$. Such level sets play a crucial role in various scientific fields, and their estimation has received significant recent interest in the fields of statistics and machine learning/pattern recognition (see below for more details). Theoretical research on this topic is mainly concerned with rates of convergence of level set estimators. While such results are interesting, they show only limited potential to be useful in practical applications. The available results do not permit statistical inference or making quantitative statements about the contour sets themselves. The contribution of this paper constitutes a significant step forward in this direction, since we establish the asymptotic normality of a class of level set estimators $\mathbb{C}_n(c)$ formed by replacing $f$ by a kernel density estimator $f_n$ in the definition of $C(c)$, in a sense that we shall soon make precise.

Received April 2008; revised September 2008.
[1]Supported in part by NSF Grant DMS-05-03908.
[2]Supported in part by NSF Grant DMS-07-06971.
*AMS 2000 subject classifications.* 60F05, 60F15, 62E20, 62G07.
*Key words and phrases.* Central limit theorem, kernel density estimator, level set estimation.

Here is our setup. Let $X_1, X_2, \ldots$ be i.i.d. with distribution function $F$ and density $f$, and consider the kernel density estimator of $f$ based on $X_1, \ldots, X_n$, $n \geq 1$,

$$f_n(x) = \frac{1}{nh_n} \sum_{i=1}^n K\left(\frac{x - X_i}{h_n^{1/d}}\right), \qquad x \in \mathbb{R}^d,$$

where $K$ is a kernel and $h_n > 0$ is a smoothing parameter. Consider the plug-in estimator

$$\mathbb{C}_n(c) = \{x : f_n(x) \geq c\}.$$

Let $G$ be a positive measure dominated by Lebesgue measure $\lambda$. Our interest is to establish the asymptotic normality of

$$\begin{aligned}(1.1) \quad d_G(\mathbb{C}_n(c), C(c)) &:= G(\mathbb{C}_n(c) \Delta C(c)) \\ &= \int_{\mathbb{R}^d} |I\{f_n(x) \geq c\} - I\{f(x) \geq c\}| \, dG(x),\end{aligned}$$

where $A \Delta B = (A \setminus B) \cup (B \setminus A)$ denotes the set-theoretic symmetric difference of two sets. Of particular interest is $G$ being the Lebesgue measure $\lambda$, as well as $G = H$ with $H$ denoting the measure having Lebesgue density $|f(x) - c|$. The latter corresponds to the so-called excess-risk which is used frequently in the classification literature, that is,

$$(1.2) \qquad d_H(\mathbb{C}_n(c), C(c)) = \int_{\mathbb{C}_n(c) \Delta C(c)} |f(x) - c| \, dx.$$

It is well known that under mild conditions $d_\lambda(\mathbb{C}_n(c), C(c)) \to 0$ in probability as $n \to \infty$, and also rates of convergence have already been derived [cf. Baíllo, Cuevas and Justel (2000), Baíllo, Cuestas-Albertos and Cuevas (2001), Cuevas, Febrero and Fraiman (2000), Baíllo (2003) and Baíllo and Cuevas (2006)]. Even more is known. Cadre (2006) derived assumptions under which for some $\mu_G > 0$ we have

$$(1.3) \qquad \sqrt{nh_n} d_G(\mathbb{C}_n(c), C(c)) \to \mu_G \qquad \text{in probability as } n \to \infty.$$

However, asymptotic normality of $d_G(\mathbb{C}_n(c), C(c))$ has not yet been considered.

Our main result says that under suitable regularity conditions there exist a normalizing sequence $\{a_{n,G}\}$ and a constant $0 < \sigma_G^2 < \infty$ such that

$$(1.4) \quad a_{n,G}\{d_G(\mathbb{C}_n(c), C(c)) - Ed_G(\mathbb{C}_n(c), C(c))\} \xrightarrow{d} \sigma_G Z \qquad \text{as } n \to \infty,$$

where $Z$ denotes a standard normal random variable. In the important special cases of $G = \lambda$ the Lebesgue measure, and $G = H$ we shall see that under suitable regularity conditions

$$(1.5) \qquad a_{n,\lambda} = \left(\frac{n}{h_n}\right)^{1/4} \quad \text{and}$$

$$(1.6) \qquad a_{n,H} = (n^3 h_n)^{1/4},$$

respectively.

In the next section we shall discuss further related work and relevant literature. In Section 2 we formulate our main result, provide some heuristics for its validity, discuss a possible statistical application and then present the proof of our result. We end Section 2 with an example and some proposals to estimate the limiting variance $\sigma_G^2$.

1.1. *Related work and literature.* Before we present our results in detail, we shall extend our overview of the literature on level set estimation to include regression level set estimation (with classification as a special case) as well as density level set estimation.

Observe that there exists a close connection between level set estimation and binary classification. The optimal (Bayes) classifier corresponds to a level set $C_\psi(0) = \{x : \psi(x) \geq 0\}$ of $\psi = pf - (1-p)g$, where $f$ and $g$ denote the Lebesgue densities of two underlying class distributions $F$ and $G$ and $p \in [0, 1]$ defines the prior probability for $f$. If an observation $X$ falls into $\{x : \psi(x) \geq 0\}$ then it is classified by the optimal classifier as coming from $F$, otherwise as coming from distribution $G$. Hall and Kang (2005) derive large sample results for this optimal classifier that are very closely related to Cadre's result (1.3). In fact, if $\text{Err}(C)$ denotes the probability of a misclassification of a binary classifier given by a set $C$, then Hall and Kong derive rates of convergence results for the quantity $\text{Err}(\hat{C}(0)) - \text{Err}(C_\psi(0))$ where $\hat{C}$ is the plug-in classifier given by $\hat{C}(0) = \{x : pf_n(x) - (1-p)g_n(x) \geq 0\}$ with $f_n$ and $g_n$ denoting the kernel estimators for $f$ and $g$, respectively. It turns out that

$$\text{Err}(\hat{C}(0)) - \text{Err}(C_\psi(0)) = \int_{\hat{C}(0) \Delta C_\psi(0)} |\psi(x)| \, dx.$$

The latter quantity is of exactly the form (1.2). The only difference is, that the function $\psi$ is not a probability density, but a (weighted) difference of two probability densities. Similarly, the plug-in estimate is a weighted difference of kernel estimates. Though the results presented here do not directly apply to this situation, the methodology used to prove them can be adapted to it in a more or less straightforward manner.

Hartigan (1975) introduced a notion of clustering via maximally connected components of density level sets. For more on this approach to clustering [see Stuetzle (2003)], and for an interesting application of this clustering approach to *astronomical sky surveys* refer to Jang (2006). Klemelä ([2004](), [2006a](), [2008]()) applies a similar point of view to develop methods for *visualizing multivariate density estimates*. Goldenshluger and Zeevi (2004) use level set estimation in the context of the Hough transform, which is a well-known *computer vision* algorithm. Certain problems in *flow cytometry* involve the statistical problem of estimating a level set for a difference of two probability densities [Roederer and Hardy (2001); see also Wand (2005)]. Further relevant applications include *detection of minefields based on arial observations*, the *analysis of seismic data*, as well as certain issues in *image segmentation*; see Huo and Lu (2004) and references therein. Another application of level set estimation is *anomaly detection* or *novelty detection*. For instance, Theiler and Cai (2003) describe how level set estimation and anomaly detection go along in the context of multispectral image analysis, where anomalous locations (pixels) correspond to unusual spectral signatures in these images. Further areas of anomaly detection include *intrusion detection* [e.g., Fan et al. (2001) and Yeung and Chow (2002)], *anomalous jet engine vibrations* [e.g., Nairac et al. (1997), Desforges, Jacob and Cooper (1998) and King et al. (2002)] or *medical imaging* [e.g., Gerig, Jomier and Chakos (2001) and Prastawa et al. (2003)] and *EEG-based seizure analysis* [Gardner et al. (2006)]. For a recent review of this area see Markou and Singh (2003).

The above list of applications of level set estimation clearly motivates the need to understand the statistical properties of level set estimators. For this reason there has been lot of recent investigation into this area. Relevant published work (not yet mentioned above) include Hartigan (1987), Polonik (1995), Cavalier (1997), Tsybakov (1997), Walther (1997), Baíllo, Cuevas and Justel (2000), Baíllo, Cuestas-Albertos and Cuevas (2001), Cuevas, Febrero and Fraiman (2000), Baíllo (2003), Tsybakov (2004), Steinwart, Hush and Scovel (2004, 2005), Gayraud and Rousseau (2005), Willett and Novak (2005, 2006), Cuevas, Gonzalez-Manteiga and Rodriguez-Casal (2006), Scott and Davenport (2006), Scott and Novak (2006), Vert and Vert (2006) and Rigollet and Vert (2008).

Finally we mention a problem closely related to that of level set estimation. This is the problem of the estimation of the support of a density, when the support is assumed to be bounded. It turns out that the methods of estimation and the techniques used to study the asymptotic properties of the estimator are very similar to those of level set estimation. Refer especially to Biau, Cadre and Pelletier (2008) and the references therein.

**2. Main result.** The rates of convergence in our main result depend on a regularity parameter $1/\gamma_g$ that describes the behavior of the slope of $g$ at the boundary set $\beta(c) = \{x \in \mathbb{R}^d : f(x) = c\}$ [see assumption (G) below]. In the important special case of $G = \lambda$ the slope of $g$ is zero, and this implies $1/\gamma_g = 0$ (or $\gamma_g = \infty$). For $G = H$ our assumptions imply that the slope of $g$ close to the boundary is bounded away from zero and infinity which says that $1/\gamma_g = 1$.

Here is our main result. *The indicated assumptions are quite technical to state and therefore for the sake of convenience they are formulated in Section 2.4 below. In particular, the integer $k \geq 1$ that appears in the statement of our theorem is defined in* (B.ii).

THEOREM 1. *Under assumptions* (D.i)–(D.ii), (K.i)–(K.ii), (G), (H) *and* (B.i)–(B.ii), *we have as* $n \to \infty$ *that*

$$a_{n,G}\{d_G(\mathbb{C}_n(c), C(c)) - Ed_G(\mathbb{C}_n(c), C(c))\} \xrightarrow{d} \sigma_G Z, \tag{2.1}$$

*where $Z$ denotes a standard normal random variable, and*

$$a_{n,G} = \left(\frac{n}{h_n}\right)^{1/4}(\sqrt{nh_n})^{1/\gamma_g}. \tag{2.2}$$

*The constant $0 < \sigma_G^2 < \infty$ is defined as in (2.57) in the case $d \geq 2$ and $k = 1$; as in (2.61) in the case $d \geq 2$ and $k \geq 2$; and as in (2.62) in the case $d = 1$ and $k \geq 2$. (The case $d = 1$ and $k = 1$ cannot occur under our assumptions.)*

REMARK 1. Write

$$\delta_n(c) = a_{n,G}\{d_G(\mathbb{C}_n(c), C(c)) - Ed_G(\mathbb{C}_n(c), C(c))\}.$$

A slight extension of the proof of our theorem shows that if $c_1, \ldots, c_m$, $m \geq 1$, are distinct positive numbers, each of which satisfies the assumptions of the theorem, then

$$(\delta_n(c_1), \ldots, \delta_n(c_m)) \xrightarrow{d} (\sigma_1 Z_1, \ldots, \sigma_m Z_m),$$

where $Z_1, \ldots, Z_m$ are independent standard normal random variables and $\sigma_1, \ldots, \sigma_m$ are as defined in the proof of the theorem.

REMARK 2. In Section 2.7 we provide an example when the variance $\sigma_G^2$ does have a closed form convenient for calculation. Such a closed form cannot be given in general, Section 2.7 also discusses some methods to estimate $\sigma_G^2$ from the data.

2.1. *Heuristics.* Before we continue with our exposition, we shall provide some heuristics to indicate why $a_n = (\frac{n}{h_n})^{1/4}$ is the correct normalizing factor in (1.5), that is, we consider the case $G = \lambda$, or $\gamma_g = \infty$. This should help the reader to understand why our theorem is true. It is well known that under certain regularity conditions we have

$$\sqrt{nh_n}\{f_n(x) - f(x)\} = O_P(1) \qquad \text{as } n \to \infty.$$

Therefore the boundary of the set $\mathbb{C}_n(c)$ can be expected to fluctuate in a band **B** with a width (roughly) of the order $O_P(\frac{1}{\sqrt{nh_n}})$ around the boundary set $\beta(c) = \{x : f(x) = c\}$. For notational simplicity we shall write $\beta = \beta(c)$. Partitioning **B** by $N = O(\frac{1}{\sqrt{nh_nh_n}}) = O(\frac{1}{\sqrt{nh_n^3}})$ regions $R_k$, $k = 1, \ldots, N$, of Lebesgue measure $\lambda(R_k) = h_n$, we can approximate $d_\lambda(\mathbb{C}_n(c), C(c))$ as

$$d_\lambda(\mathbb{C}_n(c), C(c)) \approx \sum_{k=1}^{N} \int_{R_k} |I\{f_n(x) \geq c\} - I\{f(x) \geq c\}|\, dx =: \sum_{k=1}^{N} Y_{n,k}.$$

Here we use the fact that the band **B** has width $\frac{1}{\sqrt{nh_n}}$. Writing

$$Y_{n,k} = \int_{R_k} \Lambda_n(x)\, dx$$

with

(2.3)  $$\Lambda_n(x) = |I\{f_n(x) \geq c\} - I\{f(x) \geq c\}|,$$

we see that

$$\mathrm{Var}(Y_{n,k}) = \int_{R_k}\int_{R_k} \mathrm{cov}(\Lambda_n(x), \Lambda_n(y))\, dx\, dy = O(\lambda(R_k)^2) = O(h_n^2),$$

where the $O$-terms turns out to be exact. Further, due to the nature of the kernel density estimator the variables $Y_{n,k}$ can be assumed to behave asymptotically like independent variables, since we can choose the regions $R_k$ to be disjoint. Hence, the variance of $d_\lambda(\mathbb{C}_n(c), C(c))$ can be expected to be of the order $Nh_n^2 = (\frac{h_n}{n})^{1/2}$, which motivates the normalizing factor $a_n = (\frac{n}{h_n})^{1/4}$.

2.2. *A connection to $L_p$-rates of convergence of kernel density estimates.* The following discussion on $L_p$-rates, $p \geq 1$, of convergence of kernel density estimates implicitly provides another heuristic for our result.

Consider the case $G = H_{p-1}$, where $H_{p-1}$ denotes the measure with Radon–Nikodym derivative $h_{p-1}(x) = |f(x) - c|^{p-1}$ with $p \geq 1$. Note that $H_2 = H$ with $H$ from above. Then we have the identity

(2.4)   $$\int_0^\infty H_{p-1}(\mathbb{C}_n(c)\Delta C(c))\, dc = \frac{1}{p}\int_{\mathbb{R}^d} |f_n(x) - f(x)|^p\, dx, \qquad p \geq 1.$$

The proof is straightforward [see Mason and Polonik (2008), Appendix, Detail 1]. The case $p = 1$ gives the geometrically intuitive relation

$$\int_0^\infty \lambda(\mathbb{C}_n(c)\Delta C(c))\, dc = \int_0^\infty \int_{\mathbb{C}_n(c)\Delta C(c)} dx\, dc = \int_{\mathbb{R}^d} |f_n(x) - f(x)|\, dx.$$

Assuming $f$ to be bounded, we split up the vertical axis into successive intervals $\Delta(k)$, $k = 1, \ldots, N$ of length $\approx \frac{1}{\sqrt{nh_n}}$ with midpoints $c_k$. Approximate the integral (2.4) by

$$\frac{1}{p}\int_{\mathbb{R}^d}|f_n(x)-f(x)|^p\,dx = \int_0^\infty H_{p-1}(\mathbb{C}_n(c)\Delta C(c))\,dc$$

$$\approx \sum_{k=1}^N \int_{\Delta(k)} H_{p-1}(\mathbb{C}_n(c)\Delta C(c))\,dc$$

$$\approx \frac{1}{\sqrt{nh_n}}\sum_{k=1}^N H_{p-1}(\mathbb{C}_n(c_k)\Delta C(c_k)).$$

Utilizing the $1/\sqrt{nh_n}$-rate of $f_n(x)$ we see that the last sum consists of (roughly) independent random variables. Assuming further that the variance of each (or of most) of these random variables is of the same order $a_{n,p}^{-2} = (\frac{n}{h_n})^{-1/2}(nh_n)^{-(p-1)}$ [to obtain this, apply our theorem with $\gamma_g = 1/(p-1)$] we obtain that the variance of the sum is of the order

$$\frac{a_n^{-2}}{\sqrt{nh_n}} = \left(\frac{1}{nh_n^{1-1/p}}\right)^p.$$

In other words, the normalizing factor of the $L_p$-norm of the kernel density estimator in $\mathbb{R}^d$ can be expected to be $(nh_n^{1-1/p})^{p/2} = (nh_n)^{p/2}h_n^{-1/2}$. In the case $p = 2$ this gives the normalizing factor $nh_nh_n^{-1/2} = nh_n^{1/2}$, and this coincides with the results from Rosenblatt (1975). In the special case $d = 2$ these rates can also be found in Horvath (1991).

2.3. *Possible application to online testing.* Suppose that when a certain industrial process is working properly it produces items, which may be considered as i.i.d. $\mathbb{R}^d$ random variables $X_1, X_2, \ldots$ with a known density function $f$. On the basis of a sample size $n$ taken from time to time from the production we can measure the deviation of the sample $X_1, X_2, \ldots, X_n$ from the desired distribution by looking at the discrepancy between $\lambda(C_n(c)\Delta C(c))$ and its expected value $E\lambda(C_n(c)\Delta C(c))$. The value $c$ may be chosen so that

$$P\{X \in C(c)\} = \int_{\{f(x)\geq c\}} f(x)\,dx = \alpha,$$

some typical values being $\alpha = 0.90, 0.95$ and $0.99$. We may decide to shut down the process and look for production errors if

(2.5) $\quad \sigma_\lambda^{-1}\left(\frac{n}{h_n}\right)^{1/4}|\lambda(C_n(c)\Delta C(c)) - E\lambda(C_n(c)\Delta C(c))| > 1.96.$

Otherwise as long as the estimated level set $C_n(c)$ does not deviate too much from the target level set $C(c)$ in which fraction $\alpha$ of the data should lie if the process is functioning properly, we do not disrupt production. Our central limit theorem tells us that for large enough sample sizes $n$ the probability of the event in (2.5) would be around 0.05, should, in fact, the process be working as it should. Thus using this decision rule, we make a type I error with roughly probability 0.05 if we decide to stop production, when it is actually working fine. Sometimes one might want to replace the $C_n(c)$ in the first $\lambda(C_n(c)\Delta C(c))$ in (2.5) by $C_n(c_n)$, where

$$\int_{\{f_n(x)\geq c_n\}} f_n(x)\,dx = \alpha.$$

A mechanical engineering application where this approach seems to be of some value is described in Desforges, Jacob and Cooper (1998). This application considers gearbox fault data, the collection of which is described in that paper. In fact, two classes of data were collected, corresponding to two states: a gear in good condition and a gear in bad condition, respectively. Desforges, Jacob and Cooper indicate a data analysis approach based on kernel density estimation to recognize the faulty condition. The idea is to calculate a kernel density estimator $g_m$ based on the data $X_1,\ldots,X_m$ from the gear in good condition, and then this estimator is evaluated at the data $Y_1,\ldots,Y_n$ that are sampled under a bad gear condition. Desforges, Jacob and Cooper then examine the level sets of $g_m$ in which the faulty data lie. One of their ideas is to use $\frac{1}{n}\sum_{i=1}^{n} g_m(Y_i)$, to detect the faulty condition. Their methodology is ad hoc in nature and no statistical inference procedure is proposed.

Our test procedure could be applied as follows by using $f = g_m$ (i.e., we are conditioning on $X_1,\ldots,X_m$). Set $C(c) := \{x : g_m(x) \geq c\}$, for an appropriate value of $c$, and find the corresponding set $\mathbb{C}_n(c)$ based on $Y_1,\ldots,Y_n$. Then check whether (2.5) holds. If yes, then we can conclude with at significance level 0.05 that the $Y_1,\ldots,Y_n$ stem from a different distribution than the $X_1,\ldots,X_m$. Observe that in this setup we calculate $E\lambda(\mathbb{C}_n(c)\Delta C(c))$ as well as $\sigma_\lambda^2$ by using $g_m$ as the underlying distribution. (In practice we may have to estimate these two quantities; see Section 2.7.) How this approach would work in practice is the subject of a separate paper.

2.4. *Assumptions and notation.*

ASSUMPTIONS ON THE DENSITY $f$.

(D.i) $f$ is in $C^2(\mathbb{R}^d)$ and its partial derivatives of order 1 and 2 are bounded;

(D.ii) $\inf_{x\in\mathbb{R}^d} f(x) < c < \sup_{x\in\mathbb{R}^d} f(x)$.

Notice that (D.i) implies the existence of positive constants $M$ and $A$ with

$$\sup_x f(x) \leq M < \infty \tag{2.6}$$

and

$$\frac{1}{2}\sum_{i=1}^{d}\sum_{j=1}^{d} \sup_{x\in\mathbb{R}} \left|\frac{\partial^2 f(x)}{\partial x_i\, \partial x_j}\right| =: A < \infty. \tag{2.7}$$

[Condition (D.i) implies that $f$ is uniformly continuous on $\mathbb{R}^d$ from which (2.6) follows.]

ASSUMPTIONS ON $K$.

(K.i) $K$ is a probability density function having support contained in the closed ball of radius $1/2$ centered at zero and is bounded by a constant $\kappa$.

(K.ii) $\sum_{i=1}^{d} \int_{\mathbb{R}^d} t_i K(t)\,dt = 0$.

Observe that (K.i) implies that

$$\int_{\mathbb{R}^d} |t|^2 |K(t)|\,dt = \kappa_1 < \infty. \tag{2.8}$$

ASSUMPTIONS ON THE BOUNDARY $\beta = \{x : f(x) = c\}$ FOR $d \geq 2$.

(B.i) For all $(y_1, \ldots, y_d) \in \beta$,

$$f'(y) = f'(y_1, \ldots, y_d) = \left(\frac{\partial f(y_1,\ldots,y_d)}{\partial y_1}, \ldots, \frac{\partial f(y_1,\ldots,y_d)}{\partial y_d}\right) \neq 0.$$

Define

$$I_d = \begin{cases} [0, 2\pi), & d = 2, \\ [0, \pi]^{d-2} \times [0, 2\pi), & d > 2. \end{cases}$$

The $d-1$ sphere

$$S^{d-1} = \{x \in \mathbb{R}^d : |x| = 1\}$$

can be parameterized by [e.g., see Lang (1997)]

$$x(\theta) = (x_1(\theta), \ldots, x_d(\theta)), \qquad \theta \in I_d,$$

where

$$x_1(\theta) = \cos(\theta_1),$$
$$x_2(\theta) = \sin(\theta_1)\cos(\theta_2),$$
$$x_3(\theta) = \sin(\theta_1)\sin(\theta_2)\cos(\theta_3),$$

$$\vdots$$

$$x_{d-1}(\theta) = \sin(\theta_1)\cdots\sin(\theta_{d-2})\cos(\theta_{d-1}),$$
$$x_d(\theta) = \sin(\theta_1)\cdots\sin(\theta_{d-2})\sin(\theta_{d-1}).$$

(B.ii) We assume that the boundary $\beta$ can be written as

$$\beta = \bigcup_{j=1}^{k} \beta_j \quad \text{with } \inf\{|x-y| : x \in \beta_j, y \in \beta_l\} > 0, \text{ if } j \neq l,$$

where each $\beta_j$ is diffeomorphic to $S^{d-1}$, meaning it is parameterized by a function

$$y(\theta) = (y_1(\theta), \ldots, y_d(\theta)), \qquad \theta \in I_d,$$

that is a function (depending on $j$) of the above parameterization $x(\theta)$ of $S^{d-1}$, which is 1–1 on $J_d$, the interior of $I_d$, with

$$\frac{\partial y(\theta)}{\partial \theta_i} = \left(\frac{\partial y_1(\theta)}{\partial \theta_i}, \ldots, \frac{\partial y_d(\theta)}{\partial \theta_i}\right) \neq 0, \qquad \theta \in J_d.$$

We further assume that for each $j = 1, \ldots, k$ and $i = 1, \ldots, d$, the function $\frac{\partial y(\theta)}{\partial \theta_i}$ is continuous and uniformly bounded on $J_d$.

ASSUMPTIONS ON THE BOUNDARY $\beta$ FOR $d=1$.

(B.i) $$\inf_{1 \leq i \leq k} |f'(z_i)| =: \rho_0 > 0.$$

(B.ii) $\beta = \{z_1, \ldots, z_k\}$, $k \geq 2$.

[Condition (B.i) and $f \in C^2(\mathbb{R})$ imply that the case $k = 1$ cannot occur when $d = 1$.]

ASSUMPTIONS ON $G$. (G) The measure $G$ has a bounded continuous Radon–Nikodym derivative $g$ w.r.t. Lebesgue measure $\lambda$. There exists a constant $0 < \gamma_g \leq \infty$ such that the following holds.

In the case $d \geq 2$ there exists a function $g^{(1)}(\cdot, \cdot)$ bounded on $I_d \times S^{d-1}$ such that for each $j = 1, \ldots, k$, for some $c_j \geq 0$,

$$\sup_{|z|=1} \sup_{\theta \in I_d} \left| \frac{g(y(\theta) + az)}{a^{1/\gamma_g}} - c_j g^{(1)}(\theta, z) \right| = o(1) \qquad \text{as } a \searrow 0,$$

with $0 < \sup_{|z|=1} \sup_{\theta \in I_d} |g^{(1)}(\theta, z)| < \infty$, where $y(\theta)$ is the parametrization pertaining to $\beta_j$, with at least one of the $c_j$ strictly positive.

In the case $d = 1$ there exists a function $g^{(1)}(\cdot)$ with $0 < |g^{(1)}(z_j)| < \infty$, $j = 1, \ldots, k$ such that for each $j = 1, \ldots, k$ for some $c_j \geq 0$,

$$\sup_{|z|=1} \left| \frac{g(z_j + az)}{a^{1/\gamma_g}} - c_j g^{(1)}(z_j) \right| = o(1) \qquad \text{as } a \searrow 0,$$

with at least one of the $c_j$ strictly positive. By convention, in the above statement $\frac{1}{\infty} = 0$.

ASSUMPTIONS ON $h_n$. As $n \to \infty$,

(H) $\sqrt{nh_n^{1+2/d}} \to \gamma$, with $0 \leq \gamma < \infty$ and $nh_n/\log n \to \infty$, where $\gamma = 0$ in the case $d = 1$.

**Discussion of the assumptions and some implications.**

*Discussion of assumption* (G). Measures $G$ of particular interest that satisfy assumption (G) are given by $g(x) = |f(x) - c|^p$ with $p \geq 0$, and also by $g(x) = f(x)$. The latter of course leads to the $F$-measure of the symmetric distance. The former has connections to the $L_p$-norm of the kernel density estimator (see the discussion in Section 2.2). As pointed out above in the Introduction, the choice $p = 1$ is closely connected to the excess risk from the classification literature. The choice $p = 0$ yields the Lebesgue measure of the symmetric difference.

Assumptions (B.i) and (D.i) imply that (G) holds for $|f(x) - c|^p$ with $1/\gamma_g = p$. For $g = f$ we have $1/\gamma_g = 0$ [notice that by (D.ii) we have $c > 0$].

*Discussion of smoothness assumptions on $f$.* Our smoothness assumptions on $f$ imply that $f$ has a $\gamma$-exponent with $\gamma = 1$ at the level $c$, that is, we have

$$F\{x \in \mathbb{R}^d : |f(x) - c| \leq \epsilon\} \leq C\epsilon.$$

[This fact follows from Lebesgue–Bosicovich theorem; e.g., see Cadre (2006).] This type of assumption is common in the literature of level set estimation. It was used first by Polonik (1995) in the context of density level set estimation.

*Implications of* (B) *in the case* $d \geq 2$. In the following we shall record some conventions and implications of assumption (B), which are needed in the proof of our theorem. Using the notation introduced in assumption (B), we define $\frac{\partial y(\theta)}{\partial \theta_i}$ for points on the boundary of $I_d$ to be the limit taken from points in $J_d$. In this way, we see that each vector $\frac{\partial y(\theta)}{\partial \theta_i}$ is continuous and bounded on the closure $\overline{I}_d$ of $I_d$.

Notice in the case $d \geq 2$ that for each $j = 1, \ldots, k$ and $i = 1, \ldots, d-1$,

(2.9) $$\frac{df(y(\theta))}{d\theta_i} = \frac{\partial f(y(\theta))}{\partial y_1} \frac{\partial y_1(\theta)}{\partial \theta_i} + \cdots + \frac{\partial f(y(\theta))}{\partial y_d} \frac{\partial y_d(\theta)}{\partial \theta_i} = 0,$$

where $y(\theta)$ is the parameterization pertaining to $\beta_j$. This implies that the unit vector

(2.10) $$u(\theta) = (u_1(\theta), \ldots, u_d(\theta)) := \frac{f'(y(\theta))}{|f'(y(\theta))|}$$

is normal to the tangent space of $\beta_j$ at $y(\theta)$.

From assumption (B.ii) we infer that $\beta$ is compact, which when combined with (B.i) says that

$$\inf_{(y_1,\ldots,y_d)\in\beta} |f'(y_1,\ldots,y_d)| =: \rho_0 > 0. \tag{2.11}$$

In turn, assumptions (D.ii), (B.i) and (B.ii), when combined with (2.11), imply that for each $1 \le i \le d-1$, the vector

$$\frac{\partial u(\theta)}{\partial \theta_i} = \left(\frac{\partial u_1(\theta)}{\partial \theta_i},\ldots,\frac{\partial u_d(\theta)}{\partial \theta_i}\right)$$

is uniformly bounded on $I_d$.

Consider for each $j=1,\ldots,k$, with $y(\theta)$ being the parameterization pertaining to $\beta_j$, the absolute value of the determinant,

$$\left| \det \begin{vmatrix} \frac{\partial y(\theta)}{\partial \theta_1} \\ \vdots \\ \frac{\partial y(\theta)}{\partial \theta_{d-1}} \\ u(\theta) \end{vmatrix} \right| =: \iota(\theta). \tag{2.12}$$

We can infer from (B.ii) that we have

$$\sup_{\theta \in I_d} \iota(\theta) < \infty. \tag{2.13}$$

2.5. *Proof of Theorem 1 in the case $d \ge 2$.* We shall only present a detailed proof for the case $k=1$. However, at the end we shall describe how the proof of the general $k \ge 1$ case goes. Thus for ease of notation we shall drop the subscript $j$ in the above assumptions. Also we shall assume $c_1 = 1$ in assumption (G).

We shall first show that with a suitably defined sequence of centerings $b_n$, we have

$$(n/h_n)^{1/4}(\sqrt{nh_n})^{1/\gamma_g}\{d_G(\mathbb{C}_n(c), C(c)) - b_n\} \xrightarrow{d} \sigma Z \tag{2.14}$$

for some $\sigma^2 > 0$. (For the sake of notational convenience, we write in the proof $\sigma^2 = \sigma_G^2$.) From this result we shall infer that our central limit theorem (2.1) holds. The asymptotic variance $\sigma^2$ will be defined in the course of the proof. It finally appears in (2.57) below.

Theorem 1 of Einmahl and Mason (2005) implies that when $h_n$ satisfies (H) and $f$ is bounded that for some constant $\gamma_1 > 0$

$$\limsup_{n\to\infty} \sqrt{\frac{nh_n}{\log n}} \sup_{x\in\mathbb{R}^d} |f_n(x) - Ef_n(x)| \le \gamma_1, \quad \text{a.s.} \tag{2.15}$$

It is not difficult to see that under the assumptions (D), (K) and (H) for some $\gamma_2 > 0$,

$$(2.16) \qquad \sup_{n \geq 2} \sqrt{nh_n} \sup_{x \in \mathbb{R}^d} |Ef_n(x) - f(x)| \leq \gamma_2.$$

[See Mason and Polonik (2008), Appendix, Detail 2.]

Set with $\varsigma > \sqrt{2} \vee \gamma_1$,

$$(2.17) \qquad E_n = \left\{ x : |f(x) - c| \leq \frac{\varsigma \sqrt{\log n}}{\sqrt{nh_n}} \right\}.$$

We see by (1.1), (2.15) and (2.16) that with probability 1 for all large enough $n$

$$(2.18) \qquad G(\mathbb{C}_n(c) \Delta C(c)) = \int_{E_n} |I\{f_n(x) \geq c\} - I\{f(x) \geq c\}| g(x) \, dx$$
$$=: L_n(c).$$

It turns out that rather than considering the truncated quantity $L_n(c)$ directly, it is more convenient to first study a Poissonized version of $L_n(c)$ formed by replacing $f_n(x)$ by

$$\pi_n(x) = \frac{1}{nh_n} \sum_{i=1}^{N_n} K\left(\frac{x - X_i}{h_n^{1/d}}\right),$$

where $N_n$ is a mean $n$ Poisson random variable independent of $X_1, X_2, \ldots$. [When $N_n = 0$ we set $\pi_n(x) = 0$.] Notice that

$$E\pi_n(x) = Ef_n(x).$$

We shall make repeated use of the fact following from the assumption that $K$ has support contained in the closed ball of radius $1/2$ centered at zero, that $\pi_n(x)$ and $\pi_n(y)$ are independent whenever $|x - y| > h_n^{1/d}$.

Here is the Poissonized version of $L_n(c)$ that we shall treat first. Define

$$(2.19) \qquad \Pi_n(c) = \int_{E_n} |I\{\pi_n(x) \geq c\} - I\{f(x) \geq c\}| g(x) \, dx.$$

Our goal is to infer a central limit theorem for $L_n(c)$ and thus for $G(\mathbb{C}_n(c) \Delta C(c))$ from a central limit theorem for $\Pi_n(c)$.

Set

$$(2.20) \qquad \Delta_n(x) = |I\{\pi_n(x) \geq c\} - I\{f(x) \geq c\}|.$$

The first item on this agenda is to verify that $(n/h_n)^{1/4}(\sqrt{nh_n})^{1/\gamma_g}$ is the correct sequence of norming constants. To do this we must analyze the exact asymptotic behavior of the variance of $\Pi_n(c)$. We see that

$$\operatorname{Var}(\Pi_n(c)) = \operatorname{Var}\left(\int_{E_n} \Delta_n(x)\, dG(x)\right)$$
$$= \int_{E_n}\int_{E_n} \operatorname{cov}(\Delta_n(x), \Delta_n(y))\, dG(x)\, dG(y).$$

Let

$$Y_n(x) = \left[\sum_{j \leq N_1} K\left(\frac{x - X_j}{h_n^{1/d}}\right) - EK\left(\frac{x - X}{h_n^{1/d}}\right)\right] \Big/ \sqrt{EK^2\left(\frac{x - X}{h_n^{1/d}}\right)}$$

and $Y_n^{(1)}(x), \ldots, Y_n^{(n)}(x)$ be i.i.d. $Y_n(x)$.

Clearly

$$\left(\frac{\pi_n(x) - E\pi_n(x)}{\sqrt{\operatorname{Var}(\pi_n(x))}}, \frac{\pi_n(y) - E\pi_n(y)}{\sqrt{\operatorname{Var}(\pi_n(x))}}\right)$$
$$\stackrel{d}{=} \left(\frac{\sum_{i=1}^n Y_n^{(i)}(x)}{\sqrt{n}}, \frac{\sum_{i=1}^n Y_n^{(i)}(y)}{\sqrt{n}}\right) =: (\overline{\pi}_n(x), \overline{\pi}_n(y)).$$

Set

$$c_n(x) = \frac{\sqrt{nh_n}(c - Ef_n(x))}{\sqrt{1/h_n EK^2((x-X)/h_n^{1/d})}} = \frac{\sqrt{nh_n}(c - \int_{\mathbb{R}^d} K(y)f(x - yh_n^{1/d})\, dy)}{\sqrt{1/h_n EK^2((x-X)/h_n^{1/d})}}.$$

Since $K$ has support contained in the closed ball of radius $1/2$ around zero, which implies that $\Delta_n(x)$ and $\Delta_n(y)$ are independent whenever $|x - y| > h_n^{1/d}$, we have

$$\operatorname{Var}\left(\int_{E_n} \Delta_n(x)\, dG(x)\right)$$
$$= \int_{E_n}\int_{E_n} I(|x - y| \leq h_n^{1/d}) \operatorname{cov}(\Delta_n(x), \Delta_n(y))\, dG(x)\, dG(y),$$

where now we write

$$\Delta_n(x) = \left| I\{\overline{\pi}_n(x) \geq c_n(x)\} - I\left\{0 \geq \frac{\sqrt{nh_n}(c - f(x))}{(1/h_n EK^2((x-X)/h_n^{1/d}))^{1/2}}\right\}\right|.$$

The change of variables $y = x + th_n^{1/d}$, $t \in B$, with

(2.21) $$B = \{t : |t| \leq 1\},$$

gives

(2.22) $$\operatorname{Var}\left(\int_{E_n} \Delta_n(x)\, dx\right) = h_n \int_{E_n}\int_B g_n(x, t)\, dt\, dx,$$

where

(2.23)
$$g_n(x,t) = I_{E_n}(x) I_{E_n}(x + th_n^{1/d}) \operatorname{cov}(\Delta_n(x), \Delta_n(x + th_n^{1/d})) \\ \times g(x) g(x + th_n^{1/d}).$$

For ease of notation let $a_n = a_{n,G} = (\frac{n}{h_n})^{1/4}(\sqrt{nh_n})^{1/\gamma_g}$. We intend to prove that

(2.24)
$$\lim_{n \to \infty} a_n^2 \operatorname{Var}\left(\int_{E_n} \Delta_n(x) \, dG(x)\right) \\ = \lim_{n \to \infty} a_n^2 h_n \int_{E_n} \int_B g_n(x,t) \, dt \, dx \\ = \lim_{n \to \infty} (nh_n)^{1/2 + 1/\gamma_g} \int_{E_n} \int_B g_n(x,t) \, dt \, dx \\ = \lim_{\tau \to \infty} \lim_{n \to \infty} (nh_n)^{1/2 + 1/\gamma_g} \int_{D_n(\tau)} \int_B g_n(x,t) \, dt \, dx =: \sigma^2 < \infty,$$

where

$$D_n(\tau) := \left\{ z : z = y(\theta) + \frac{su(\theta)}{\sqrt{nh_n}}, \theta \in I_d, |s| \leq \tau \right\}.$$

The set $D_n(\tau)$ forms a band around the surface $\beta$ of thickness $\frac{2\tau}{\sqrt{nh_n}}$.

Recall the definition of $B$ in (2.21). Since $\beta$ is a closed submanifold of $\mathbb{R}^d$ without boundary the *tubular neighborhood theorem* [see Theorem 11.4 on page 93 of Bredon (1993)] says that for all $\delta > 0$ sufficiently small for each $x \in \beta + \delta B$ there is an unique $\theta \in I_d$ and $|s| \leq \delta$ such that $x = y(\theta) + su(\theta)$. This, in turn, implies that for all $\delta > 0$ sufficiently small

$$\{y(\theta) + su(\theta) : \theta \in I_d \text{ and } |s| \leq \delta\} = \beta + \delta B.$$

In particular, we see by using (H) that for all large enough $n$

(2.25) $$D_n(\tau) = \beta + \frac{\tau}{\sqrt{nh_n}} B, \qquad \text{where } B = \{z : |z| \leq 1\}.$$

Moreover, it says that $x = y(\theta) + su(\theta)$, $\theta \in I_d$ and $|s| \leq \delta$ is a well-defined parameterization of $\beta + \delta B$, and it validates the change of variables in the integrals below.

We now turn to the proof of (2.24). Let

$$\rho_n(x, x + th_n^{1/d}) = \operatorname{Cov}(\overline{\pi}_n(x), \overline{\pi}_n(x + th_n^{1/d})) \\ = \frac{h_n^{-1} E[K((x-X)/h_n^{1d}) K((x-X)/h_n^{1d} + t)]}{\sqrt{h_n^{-1} E K^2((x-X)/h_n^{1d}) h_n^{-1} E K^2((x-X)/h_n^{1d} + t)}}.$$

It is routine to show that for each $\theta \in I_d$, $|s| \leq \tau$, $x = y(\theta) + \frac{su(\theta)}{\sqrt{nh_n}}$ and $t \in B$ we have as $n \to \infty$ that

$$\rho_n(x, x + th_n^{1/d}) = \rho_n\left(y(\theta) + \frac{su(\theta)}{\sqrt{nh_n}}, y(\theta) + \frac{su(\theta)}{\sqrt{nh_n}} + th_n^{1/d}\right) \to \rho(t),$$

where

$$\rho(t) := \frac{\int_{\mathbb{R}^d} K(u) K(u+t)\, du}{\int_{\mathbb{R}^d} K^2(u)\, du}.$$

[See Mason and Polonik (2008), Appendix, Detail 4.] Notice that $\rho(t) = \rho(-t)$. One can then infer by the central limit theorem that for each $\theta \in I_d$, $|s| \leq \tau$ and $t \in B$,

$$(\overline{\pi}_n(x), \overline{\pi}_n(x + th_n^{1/d}))$$
(2.26)
$$= \left(\overline{\pi}_n\left(y(\theta) + \frac{su(\theta)}{\sqrt{nh_n}}\right), \overline{\pi}_n\left(y(\theta) + \frac{su(\theta)}{\sqrt{nh_n}} + th_n^{1/d}\right)\right)$$
$$\xrightarrow{d} (Z_1, \rho(t)Z_1 + \sqrt{1 - \rho^2(t)}Z_2),$$

where $Z_1$ and $Z_2$ are independent standard normal random variables.

We also get by using our assumptions and straightforward Taylor expansions that for $|s| \leq \tau$, $u = su(\theta)$, $x = y(\theta) + \frac{u}{\sqrt{nh_n}}$ and $\theta \in I_d$

(2.27)
$$\begin{aligned}
c_n(x) &= c_n\left(y(\theta) + \frac{u}{\sqrt{nh_n}}\right) \\
&= \frac{\sqrt{nh_n}(c - Ef_n(y(\theta) + u/\sqrt{nh_n}))}{\sqrt{1/h_n EK^2((y(\theta) + u/\sqrt{nh_n} - X)/h_n^{1/d})}} \\
&\xrightarrow{n\to\infty} -\frac{f'(y(\theta)) \cdot u}{\sqrt{f(y(\theta))}\|K\|_2} \\
&= -\frac{s|f'(y(\theta))|}{\sqrt{c}\|K\|_2} =: c(s, \theta, 0)
\end{aligned}$$

and similarly since $\sqrt{nh_n^{1+2/d}} \to \gamma$,

$$\begin{aligned}
c_n(x + th_n^{1/d}) &\xrightarrow{n\to\infty} -\frac{f'(y(\theta)) \cdot (u+t)}{\sqrt{f(y(\theta))}\|K\|_2} \\
&= -\frac{s|f'(y(\theta))|}{\sqrt{c}\|K\|_2} - \frac{\gamma f'(y(\theta)) \cdot t}{\sqrt{c}\|K\|_2} \\
&= : c(s, \theta, \gamma t).
\end{aligned}$$

We also have

$$\frac{\sqrt{nh_n}(c - f(x))}{\sqrt{1/h_n EK^2((x-X)/h_n^{1/d})}} \overset{n \to \infty}{\to} c(s, \theta, 0)$$

and

$$\frac{\sqrt{nh_n}(c - f(x + th_n^{1/d}))}{\sqrt{1/h_n EK^2((x-X)/h_n^{1/d})}} \overset{n \to \infty}{\to} c(s, \theta, \gamma t).$$

[See Mason and Polonik (2008), Appendix, Detail 5.] Hence by (2.26) and (G) for $y(\theta) \in \beta$,

(2.28)
$$\begin{aligned}
&(nh_n)^{1/\gamma_g} g_n(x, t) \\
&= (nh_n)^{1/\gamma_g} g_n\left(y(\theta) + \frac{u}{\sqrt{nh_n}}, t\right) \\
&= I_{E_n}\left(y(\theta) + \frac{u}{\sqrt{nh_n}}\right) I_{E_n}\left(y(\theta) + \frac{u}{\sqrt{nh_n}} + th_n^{1/d}\right) \\
&\quad \times \operatorname{cov}\left(\Delta_n\left(y(\theta) + \frac{u}{\sqrt{nh_n}}\right), \Delta_n\left(y(\theta) + \frac{u}{\sqrt{nh_n}} + th_n^{1/d}\right)\right) \\
&\quad \times (nh_n)^{1/\gamma_g} g\left(y(\theta) + \frac{u}{\sqrt{nh_n}}\right) g\left(y(\theta) + \frac{u}{\sqrt{nh_n}} + th_n^{1/d}\right) \\
&\overset{n \to \infty}{\to} \operatorname{cov}(|I\{Z_1 \geq c(s, \theta, 0)\} - I\{0 \geq c(s, \theta, 0)\}|, \\
&\qquad\qquad |I\{\rho(t) Z_1 + \sqrt{1 - \rho^2(t)} Z_2 \geq c(s, \theta, \gamma t)\} \\
&\qquad\qquad\qquad - I\{0 \geq c(s, \theta, \gamma t)\}|) \\
&\quad \times |s|^{1/\gamma_g} g^{(1)}(\theta, u(\theta)) |su(\theta) + \gamma t|^{1/\gamma_g} g^{(1)}\left(\theta, \frac{su(\theta) + \gamma t}{|su(\theta) + \gamma t|}\right) \\
&=: \Gamma(\theta, s, t).
\end{aligned}$$

Using the change of variables

(2.29) $$x_1 = y_1(\theta) + \frac{su_1(\theta)}{\sqrt{nh_n}}, \ldots, x_d = y_d(\theta) + \frac{su_d(\theta)}{\sqrt{nh_n}},$$

we get

$$\int_{D_n(\tau)} \int_B g_n(x, t) \, dt \, dx = \int_{-\tau}^{\tau} \int_{I_d} \int_B g_n\left(y(\theta) + \frac{su(\theta)}{\sqrt{nh_n}}, t\right) |J_n(\theta, s)| \, dt \, d\theta \, ds,$$

where

$$|J_n(\theta, s)| = \left| \det \begin{Vmatrix} \frac{\partial y(\theta)}{\partial \theta_1} + \frac{1}{\sqrt{nh_n}} \frac{\partial u(\theta)}{\partial \theta_1} \\ \vdots \\ \frac{\partial y(\theta)}{\partial \theta_{d-1}} + \frac{1}{\sqrt{nh_n}} \frac{\partial u(\theta)}{\partial \theta_{d-1}} \\ \frac{u(\theta)}{\sqrt{nh_n}} \end{Vmatrix} \right|. \quad (2.30)$$

Clearly, with $\iota(\theta)$ as in (2.12),

$$\sqrt{nh_n}|J_n(\theta, s)| \to \iota(\theta). \quad (2.31)$$

Under our assumptions we have $\sqrt{nh_n}|J_n(\theta, s)|$ is uniformly bounded in $n \geq 1$ and $(\theta, s) \in I_d \times [-\tau, \tau]$. Also by using (G) we see that for all $n$ large enough $(nh_n)^{1/\gamma_g} g_n$ is bounded on $I_d \times B$. Thus since $(nh_n)^{1/\gamma_g} g_n$ and $\sqrt{nh_n}|J_n|$ are eventually bounded on the appropriate domains, and (2.28) and (2.31) hold, we get by the dominated convergence theorem and (G) that

$$(nh_n)^{1/2+1/\gamma_g} \int_{D_n(\tau)} \int_B g_n(x, t) \, dt \, dx$$

$$= (nh_n)^{1/2+1/\gamma_g} \int_{-\tau}^{\tau} \int_{I_d} \int_B g_n\left(y(\theta) + \frac{su(\theta)}{\sqrt{nh_n}}, t\right) |J_n(\theta, s)| \, dt \, d\theta \, ds \quad (2.32)$$

$$\to \int_{-\tau}^{\tau} \int_{I_d} \int_B \Gamma(\theta, s, t) \iota(\theta) \, dt \, d\theta \, ds \quad \text{as } n \to \infty.$$

We claim that as $\tau \to \infty$ we have

$$\int_{-\tau}^{\tau} \int_{I_d} \int_B \Gamma(\theta, s, t) \iota(\theta) \, dt \, d\theta \, ds$$

$$\to \int_{-\infty}^{\infty} \int_{I_d} \int_B \Gamma(\theta, s, t) \iota(\theta) \, dt \, d\theta \, ds =: \sigma^2 < \infty \quad (2.33)$$

and

$$\lim_{\tau \to \infty} \limsup_{n \to \infty} (nh_n)^{1/2+1/\gamma_g} \int_{D_n^C(\tau) \cap E_n} \int_B g_n(x, t) \, dt \, dx = 0, \quad (2.34)$$

which in light of (2.32) implies that the limit in (2.24) is equal to $\sigma^2$ as defined in (2.33).

First we show (2.33). Consider

$$\Gamma^+(\tau) := \int_0^{\tau} \int_{I_d} \int_B \Gamma(\theta, s, t) \iota(\theta) \, dt \, d\theta \, ds.$$

We shall show existence and finiteness of the limit $\lim_{\tau\to\infty}\Gamma^+(\tau)$. Similar arguments apply to

$$\lim_{\tau\to\infty}\Gamma^-(\tau) := \lim_{\tau\to\infty}\int_{-\tau}^0\int_{I_d}\int_B \Gamma(\theta,s,t)\iota(\theta)\,dt\,d\theta\,ds < \infty.$$

Observe that when $s \geq 0$,

$$|I\{Z_1 \geq c(s,\theta,0)\} - I\{0 \geq c(s,\theta,0)\}| = I\{Z_1 < c(s,\theta,0)\}$$

and with $\Phi$ denoting the cdf of a standard normal distribution we write

$$E(I\{Z_1 < c(s,\theta,0)\}) = \Phi(c(s,\theta,0)).$$

Hence by taking into account (2.28), the assumed finiteness of $\sup_{|z|=1}\sup_\theta g^{(1)}(\theta,z)$, and using the elementary inequality

$$|\operatorname{cov}(X,Y)| \leq 2E|X|, \qquad \text{whenever } |Y| \leq 1,$$

we get for all $s \geq 0$ and some $c_1 > 0$ that

$$(2.35) \qquad |\Gamma(\theta,s,t)| \leq c_1|s|^{1/\gamma_g}(|s|^{1/\gamma_g} + \gamma^{1/\gamma_g})\Phi(c(s,\theta,0)).$$

The lower bound (2.11) implies the existence of a constant $\tilde{c} > 0$ such that

$$\Phi(c(s,\theta,0)) = \Phi\left(-\frac{s|f'(y(\theta))|}{\sqrt{c}\|K\|_2}\right) \leq \Phi(-\tilde{c}s).$$

Together with (2.35) and (2.13) it follows that for some $\overline{c} > 0$ we have

$$\lim_{\tau\to\infty}\Gamma^+(\tau) \leq \overline{c}\lim_{\tau\to\infty}\int_0^\tau |s|^{1/\gamma_g}(|s|^{1/\gamma_g} + \gamma^{1/\gamma_g})(1-\Phi(\tilde{c}s))\,ds < \infty.$$

Similarly,

$$|\Gamma^+(\infty) - \Gamma^+(\tau)| \leq \overline{c}\int_\tau^\infty |s|^{1/\gamma_g}(|s|^{1/\gamma_g} + \gamma^{1/\gamma_g})(1-\Phi(\tilde{c}s))\,ds \to 0$$

as $\tau \to \infty$.

This validates claim (2.33).

Next we turn to the proof of (2.34). Recall the definition of $g_n(x,t)$ in (2.23). Notice that for all $n$ large enough, we have

$$(nh_n)^{1/2+1/\gamma_g}\int_{D_n^C(\tau)\cap E_n}\int_B g_n(x,t)\,dt\,dx$$

$$\leq \sqrt{nh_n}\int_{D_n^C(\tau)\cap E_n}\int_B |\operatorname{cov}(\Delta_n(x),\Delta_n(x+th_n^{1/d}))|$$

$$(2.36) \qquad\qquad \times (nh_n)^{1/\gamma_g}g(x)g(x+th_n^{1/d})\,dt\,dx$$

$$\leq \sqrt{nh_n}\int_{D_n^C(\tau)\cap E_n}\int_B (\operatorname{Var}(\Delta_n(x)))^{1/2}$$

$$(2.37) \qquad\qquad \times (nh_n)^{1/\gamma_g}g(x)g(x+th_n^{1/d})\,dt\,dx.$$

The last inequality uses the fact that $\Delta_n(x+th_n^{1/d}) \leq 1$ and thus $\text{Var}(\Delta_n(x+th_n^{1/d})) \leq 1$. Applying the inequality

(2.38) $$|I\{a \geq b\} - I\{0 \geq b\}| \leq I\{|a| \geq |b|\},$$

we obtain that

$$\begin{aligned}\text{Var}(\Delta_n(x)) &= \text{Var}(|I\{\pi_n(x) \geq c\} - I\{f(x) \geq c\}|) \\ &\leq E(I\{|\pi_n(x) - f(x)| \geq |c - f(x)|\})^2 \\ &= P\{|\pi_n(x) - f(x)| \geq |c - f(x)|\}.\end{aligned}$$

Thus we get that

$$\sqrt{nh_n} \int_{D_n^C(\tau) \cap E_n} \int_B \sqrt{\text{Var}(\Delta_n(x))}(nh_n)^{1/\gamma_g} g(x)g(x+th_n^{1/d})\,dt\,dx$$

(2.39) $$\leq \sqrt{nh_n} \int_{D_n^C(\tau) \cap E_n} \int_B \sqrt{P\{|\pi_n(x) - f(x)| \geq |c - f(x)|\}}$$

$$\times (nh_n)^{1/\gamma_g} g(x)g(x+th_n^{1/d})\,dt\,dx.$$

We must bound the probability inside the integral. For this purpose we need a lemma.

LEMMA 2.1. *Let $Y, Y_1, Y_2, \ldots$ be i.i.d. with mean $\mu$ and bounded by $0 < M < \infty$. Independent of $Y_1, Y_2, \ldots$ let $N_n$ be a Poisson random variable with mean $n$. For any $v \geq 2(e30)^2 EY^2$ and with $d = e30M$ we have for all $\lambda > 0$,*

(2.40) $$P\left\{\sum_{i=1}^{N_n} Y_i - n\mu \geq \lambda\right\} \leq \exp\left(-\frac{\lambda^2/2}{nv + d\lambda}\right).$$

PROOF. Let $N$ be a Poisson random variable with mean 1 independent of $Y_1, Y_2, \ldots$ and let

$$\omega = \sum_{i=1}^N Y_i.$$

Clearly if $\omega_1, \ldots, \omega_n$ are i.i.d. $\omega$, then

$$\sum_{i=1}^{N_n} Y_i - n\mu \stackrel{d}{=} \sum_{i=1}^n (\omega_i - \mu).$$

Our aim is to use Bernstein's inequality to prove (2.40). Notice that for any integer $r \geq 2$,

(2.41) $$E|\omega - \mu|^r = E\left|\sum_{i=1}^N Y_i - \mu\right|^r.$$

At this point we need the following fact, which is Lemma 2.3 of Giné, Mason and Zaitsev (2003).

FACT 1. If, for each $n \geq 1$, $\zeta, \zeta_1, \zeta_2, \ldots, \zeta_n, \ldots$, are independent identically distributed random variables, $\zeta_0 = 0$, and $\eta$ is a Poisson random variable with mean $\gamma > 0$ and independent of the variables $\{\zeta_i\}_{i=1}^\infty$, then, for every $p \geq 2$,

$$(2.42) \qquad E\left|\sum_{i=0}^\eta \zeta_i - \gamma E\zeta\right|^p \leq \left(\frac{15p}{\log p}\right)^p \max[(\gamma E\zeta^2)^{p/2}, \gamma E|\zeta|^p].$$

Applying inequality (2.42) to (2.41) gives for $r \geq 2$

$$E|\omega - \mu|^r = E\left|\sum_{i=1}^N Y_i - \mu\right|^r \leq \left(\frac{15r}{\log r}\right)^r \max[(EY^2)^{r/2}, E|Y|^r].$$

Now

$$\max[(EY^2)^{r/2}, E|Y|^r] \leq \max[(EY^2)(EY^2)^{r/2-1}, (EY^2)M^{r-2}]$$
$$\leq EY^2 M^{r-2}.$$

Moreover, since $\log 2 \geq 1/2$, we get

$$E|\omega - \mu|^r \leq (30r)^r EY^2 M^{r-2}.$$

By Stirling's formula [see page 864 of Shorack and Wellner (1986)]

$$r^r \leq e^r r!.$$

Thus

$$E|\omega - \mu|^r \leq (e30r)^r EY^2 M^{r-2} \leq \frac{2(e30)^2 EY^2}{2} r!(e30M)^{r-2} \leq \frac{v}{2} r! d^{r-2},$$

where $v \geq 2(e30)^2 EY^2$ and $d = e30M$. Thus by Bernstein's inequality [see page 855 of Shorack and Wellner (1986)] we get (2.40).

Here is how Lemma 2.1 is used. Let $Y_i = K(\frac{x-X_i}{h_n^{1/d}})$. Since by assumption both $K$ and $f$ are bounded, and $K$ has support contained in the closed ball of radius $1/2$ around zero, we obtain that for some $D_0 > 0$ and all $n \geq 1$,

$$\sup_{x \in \mathbb{R}^d} E\left[K\left(\frac{x-X}{h_n^{1/d}}\right)\right]^2 \leq D_0 h_n.$$

Consider $z \geq a/\sqrt{nh_n}$ for some $a > 0$. With this choice, and since $\sup_x |Ef_n(x) - f(x)| \leq A_1 h^{2/d} \leq \frac{a}{2\sqrt{nh_n}}$ for $n$ large enough by using (H) [see Mason and Polonik

(2008), Appendix, Detail 2], we have

$$P\{\pi_n(x) - f(x) \geq z\} = P\{\pi_n(x) - Ef_n(x) \geq z - (Ef_n(x) - f(x))\}$$
$$\leq P\left\{\pi_n(x) - Ef_n(x) \geq z - \frac{1}{2}\frac{a}{\sqrt{nh_n}}\right\}$$
$$\leq P\left\{\pi_n(x) - Ef_n(x) \geq \frac{z}{2}\right\}$$

for $z \geq a/\sqrt{nh_n}$ and $n$ large enough. We get then from inequality (2.40) that for $n \geq 1$, all $z > 0$ that for some constants $D_1$ and $D_2$

$$P\left\{\pi_n(x) - Ef_n(x) \geq \frac{z}{2}\right\} = P\left\{\sum_{i=1}^{N_n} K\left(\frac{x - X_i}{h_n^{1/d}}\right) - nEK\left(\frac{x - X}{h_n^{1/d}}\right) \geq \frac{nh_n z}{2}\right\}$$
$$\leq \exp\left(-\frac{(nh_n)^2 z^2}{D_1 nh_n + D_2 nh_n z}\right) = \exp\left(-\frac{nh_n z^2}{D_1 + D_2 z}\right).$$

We see that for some $a > 0$ for all $z \geq a/\sqrt{nh_n}$ and $n$ large enough,

$$\frac{nh_n z^2}{D_1 + D_2 z} \geq \sqrt{nh_n} z.$$

Observe that for $0 \leq z \leq a/\sqrt{nh_n}$,

$$\exp(a)\exp(-\sqrt{nh_n} z) \geq \exp(a)\exp(-a) = 1 \geq P\{\pi_n(x) - f(x) \geq z\}.$$

Therefore by setting $A = \exp(a)$ we get for all large enough $n \geq 1$, $z > 0$ and $x$,

$$P\{\pi_n(x) - f(x) \geq z\} \leq A\exp(-\sqrt{nh_n} z).$$

In the same way, for all large enough $n \geq 1$, $z > 0$ and $x$,

$$P\{\pi_n(x) - f(x) \leq -z\} \leq A\exp(-\sqrt{nh_n} z).$$

Notice these inequalities imply that for all large enough $n \geq 1$, $z > 0$ and $x$,

$$(2.43) \quad \sqrt{P\{|\pi_n(x) - f(x)| \geq |c - f(x)|\}} \leq \sqrt{A}\exp\left(-\frac{\sqrt{nh_n}|c - f(x)|}{2}\right).$$

Returning to the proof of (2.34), from (2.36), (2.37), (2.39) and (2.43) we get that for all large enough $n \geq 1$,

$$(nh_n)^{1/2 + 1/\gamma_g} \int_{D_n^C(\tau) \cap E_n} \int_B g_n(x, t)\, dt\, dx$$
$$\leq \sqrt{nh_n}\sqrt{A} \int_{D_n^C(\tau) \cap E_n} \int_B e^{-\sqrt{nh_n}|c - f(x)|/2} (nh_n)^{1/\gamma_g} g(x) g(x + th_n^{1/d})\, dt\, dx,$$

which equals

$$\lambda(B)\sqrt{nh_n}\int_{D_n^C(\tau)\cap E_n}\varphi_n(x)\,dx,$$

where

$$\varphi_n(x) = \sqrt{A}\exp\left(-\frac{\sqrt{nh_n}|c-f(x)|}{2}\right)(nh_n)^{1/\gamma_g}g(x)g(x+th_n^{1/d}).$$

Our assumptions imply that for some $0 < \eta < 1$ for all $1 \leq |s| \leq \varsigma\sqrt{\log n}$ and $n$ large

$$\frac{\sqrt{nh_n}}{2}\left|c - f\left(y(\theta) + \frac{su(\theta)}{\sqrt{nh_n}}\right)\right| \geq \eta|s|.$$

[See Mason and Polonik (2008), Appendix, Detail 3.] We get using the change of variables (2.29) that for all $\tau > 1$,

$$\int_{D_n^C(\tau)\cap E_n}\varphi_n(x)\,dx = \int_{\tau \leq |s| \leq \varsigma\sqrt{\log n}}\int_{I_d}\varphi_n\left(y(\theta) + \frac{su(\theta)}{\sqrt{nh_n}}\right)|J_n(\theta,s)|\,d\theta\,ds.$$

Thus, by our assumptions [refer to the remarks after (2.31) and assumption (G)] there exists a constant $C > 0$, such that for all large enough $\tau$ and $n$

$$\int_{D_n^C(\tau)\cap E_n}\varphi_n(x)\,dx$$
$$\leq \frac{C}{\sqrt{nh_n}}$$
$$\times \int_{\tau \leq |s| \leq \varsigma\sqrt{\log n}}\int_{I_d}|s|^{2/\gamma_g}$$
$$\times \exp\left(-\frac{\sqrt{nh_n}|c - f(y(\theta) + su(\theta)/\sqrt{nh_n})|}{2}\right)d\theta\,ds$$
$$\leq \frac{C}{\sqrt{nh_n}}\int_{\tau \leq |s| \leq \varsigma\sqrt{\log n}}\int_{I_d}\exp\left(-\frac{\eta|s|}{2}\right)d\theta\,ds.$$

Thus

(2.44) $$\int_{D_n^C(\tau)\cap E_n}\int_B\varphi_n(x)\,dx \leq \frac{4\pi^{d-1}C\exp(-\eta\tau/2)}{\eta\sqrt{nh_n}}.$$

Therefore after inserting all of the above bounds we get that

$$(nh_n)^{1/2+1/\gamma_g}\int_{D_n^C(\tau)\cap E_n}\int_B g_n(x,t)\,dx\,dt \leq \frac{4\pi^{d-1}C\exp(-\eta\tau/2)}{\eta}$$

and hence we readily conclude that (2.34) holds.

Putting everything together we get that as $n \to \infty$,

$$a_n^2 \operatorname{Var}(\Pi_n(c)) \to \sigma^2 \tag{2.45}$$

with $\sigma^2$ defined as in (2.33). For future use, we point out that we can infer by (2.24), (2.34) and (2.45) that for all $\varepsilon > 0$ there exist a $\tau_0$ and an $n_0 \geq 1$ such that for all $\tau \geq \tau_0$ and $n \geq n_0$

$$|\sigma_n^2(\tau) - \sigma^2| < \varepsilon, \tag{2.46}$$

where

$$\sigma_n^2(\tau) = \operatorname{Var}\left(\left(\frac{n}{h_n}\right)^{1/4} (\sqrt{nh_n})^{1/\gamma_g} \int_{D_n(\tau)} \Delta_n(x) g(x)\, dx\right). \tag{2.47}$$

Our next goal is to de-Poissonize by applying the following version of a theorem in Beirlant and Mason (1995).

LEMMA 2.2. *Let $N_{1,n}$ and $N_{2,n}$ be independent Poisson random variables with $N_{1,n}$ being $\operatorname{Poisson}(n\beta_n)$ and $N_{2,n}$ being $\operatorname{Poisson}(n(1-\beta_n))$ where $\beta_n \in (0,1)$. Denote $N_n = N_{1,n} + N_{2,n}$ and set*

$$U_n = \frac{N_{1,n} - n\beta_n}{\sqrt{n}} \quad and \quad V_n = \frac{N_{2,n} - n(1-\beta_n)}{\sqrt{n}}.$$

*Let $\{S_n\}_{n=1}^{\infty}$ be a sequence of random variables such that:*

(i) *for each $n \geq 1$, the random vector $(S_n, U_n)$ is independent of $V_n$,*
(ii) *for some $\sigma^2 < \infty$, $S_n \xrightarrow{d} \sigma Z$, as $n \to \infty$,*
(iii) *$\beta_n \to 0$, as $n \to \infty$.*

*Then, for all $x$,*

$$P\{S_n \leq x \mid N_n = n\} \to P\{\sigma Z \leq x\}.$$

The proof follows along the same lines as Lemma 2.4 in Beirlant and Mason (1995). [See Mason and Polonik (2008), Appendix, Detail 6.]

We shall now use this de-Poissonization lemma to complete the proof of our theorem. Recall the definitions of $L_n(c)$ and $\Pi_n(c)$ in (2.18) and (2.19), respectively. Noting that $D_n(\tau) \subset E_n$ for all large enough $n \geq 1$, we see that

$$a_n(L_n(c) - E\Pi_n(c))$$
$$= a_n \int_{D_n(\tau)} \{|I\{f_n(x) \geq c\} - I\{f(x) \geq c\}| - E\Delta_n(x)\} g(x)\, dx$$
$$\quad + a_n \int_{D_n(\tau)^C \cap E_n} \{|I\{f_n(x) \geq c\} - I\{f(x) \geq c\}| - E\Delta_n(x)\} g(x)\, dx$$
$$=: T_n(\tau) + R_n(\tau).$$

We can control the $R_n(\tau)$ piece of this sum using the inequality, which follows from Lemma 2.3 below,

$$E(R_n(\tau))^2$$
$$(2.48) \quad \leq 2a_n^2 \operatorname{Var}\left(\int_{D_n(\tau)^C \cap E_n} |I\{\pi_n(x) \geq c\} - I\{f(x) \geq c\}| g(x)\, dx\right)$$
$$= 2(nh_n)^{1/2+1/\gamma_g} \int_{D_n^C(\tau) \cap E_n} \int_B g_n(x,t)\, dt\, dx,$$

which goes to zero as $n \to \infty$ and $\tau \to \infty$ as we proved in (2.34).

The needed inequality is a special case of the following result in Giné, Mason and Zaitsev (2003). We say that a set $D$ is a (commutative) semigroup if it has a commutative and associative operation, in our case sum, with a zero element. If $D$ is equipped with a $\sigma$-algebra $\mathcal{D}$ for which the sum, $+: (D \times D, \mathcal{D} \otimes \mathcal{D}) \mapsto (D, \mathcal{D})$, is measurable, then we say the $(D, \mathcal{D})$ is a measurable semigroup.

LEMMA 2.3. *Let $(D, \mathcal{D})$ be a measurable semigroup; let $Y_0 = 0 \in D$ and let $Y_i$, $i \in \mathbf{N}$, be independent identically distributed $D$-valued random variables; for any given $n \in \mathbf{N}$, let $\eta$ be a Poisson random variable with mean $n$ independent of the sequence $\{Y_i\}$; and let $B \in \mathcal{D}$ be such that $P\{Y_1 \in B\} \leq 1/2$. If $G: D \mapsto \mathbf{R}$ is nonnegative and $\mathcal{D}$-measurable, then*

$$(2.49) \quad EG\left(\sum_{i=0}^n I(Y_i \in B)Y_i\right) \leq 2EG\left(\sum_{i=0}^\eta I(Y_i \in B)Y_i\right).$$

Next we consider $T_n(\tau)$. Observe that

$$(2.50) \quad (S_n(\tau)|N_n = n) \stackrel{d}{=} \frac{T_n(\tau)}{\sigma_n(\tau)},$$

where as above $N_n$ denotes a Poisson random variable with mean $n$,

$$S_n(\tau) = \frac{a_n \int_{D_n(\tau)} \{\Delta_n(x) - E\Delta_n(x)\} g(x)\, dx}{\sigma_n(\tau)},$$

and $\sigma_n^2(\tau)$ is defined as in (2.47). We shall apply Lemma 2.2 to $S_n(\tau)$ with

$$N_{1,n} = \sum_{i=1}^{N_n} 1\{X_i \in D_n(\tau + \sqrt{n}h_n)\},$$

$$N_{2,n} = \sum_{i=1}^{N_n} 1\{X_i \notin D_n(\tau + \sqrt{n}h_n)\}$$

and
$$\beta_n = P\{X_i \in D_n(\tau + \sqrt{n}h_n)\}.$$

We first need to verify that as $n \to \infty$
$$S_n(\tau) = \frac{a_n \int_{D_n(\tau)} \{\Delta_n(x) - E\Delta_n(x)\} g(x)\, dx}{\sigma_n(\tau)} \xrightarrow{d} Z.$$

To show this we require the following special case of Theorem 1 of Shergin (1990).

FACT 2 [Shergin (1990)]. Let $\{X_{\mathbf{i},n} : \mathbf{i} \in \mathbf{Z}^d\}$ denote a triangular array of mean zero $m$-dependent random fields, and let $\mathcal{J}_n \subset \mathbf{Z}^d$ be such that:

(i) $\mathrm{Var}(\sum_{\mathbf{i} \in \mathcal{J}_n} X_{\mathbf{i},n}) \to 1$, as $n \to \infty$, and
(ii) for some $2 < s < 3$, $\sum_{\mathbf{i} \in \mathcal{J}_n} E|X_{\mathbf{i},n}|^s \to 0$, as $n \to \infty$.

Then
$$\sum_{\mathbf{i} \in \mathcal{J}_n} X_{\mathbf{i},n} \xrightarrow{d} Z,$$

where $Z$ is a standard normal random variable.

We use Shergin's result as follows. Under our regularity conditions, for each $\tau > 0$ there exist positive constants $d_1, \ldots, d_5$ such that for all large enough $n$,

(2.51) $$|D_n(\tau)| \leq \frac{d_1}{\sqrt{nh_n}};$$

(2.52) $$d_2 \leq \sigma_n(\tau) \leq d_3.$$

Clearly (2.52) follows from (2.46), and it is not difficult to see (2.51). For details see Mason and Polonik (2008), Appendix, Detail 7. There it is also shown that for each such integer $n \geq 1$ there exists a partition $\{R_{\mathbf{i}}, \mathbf{i} \in \mathcal{J}_n \subset Z^d\}$ of $D_n(\tau)$ such that for each $\mathbf{i} \in \mathcal{J}_n$

(2.53) $$|R_{\mathbf{i}}| \leq d_4 h_n,$$

where

(2.54) $$|\mathcal{J}_n| =: m_n \leq \frac{d_5}{\sqrt{nh_n^3}}.$$

Define
$$X_{\mathbf{i},n} = \frac{a_n \int_{R_{\mathbf{i}}} \{\Delta_n(x) - E\Delta_n(x)\} g(x)\, dx}{\sigma_n(\tau)}, \qquad \mathbf{i} \in \mathcal{J}_n.$$

It is straightforward to see that $X_{\mathbf{i},n}$ can be extended to a 1-dependent random field on $\mathbf{Z}^d$. [See Mason and Polonik (2008), Appendix, Detail 7.]

Notice that by (G) there exists a constant $A > 0$ such that for all $x \in D_n(\tau)$,
$$|g(x)| \leq A(\sqrt{nh_n})^{-1/\gamma_g}.$$

Recalling that $a_n = a_{n,G} = (\frac{n}{h_n})^{1/4}(\sqrt{nh_n})^{1/\gamma_g}$ we thus obtain for all for $\mathbf{i} \in \mathcal{J}_n$,
$$|X_{\mathbf{i},n}| \leq \frac{a_n 2A|R_{\mathbf{i}}|A(\sqrt{nh_n})^{-1/\gamma_g}}{\sigma_n(\tau)}$$
$$\leq \frac{2d_4 A h_n}{d_2}\left(\frac{n}{h_n}\right)^{1/4} = \frac{2Ad_4}{d_2}(nh_n^3)^{1/4}.$$

Therefore,
$$\sum_{\mathbf{i} \in \mathcal{J}_n} E|X_{\mathbf{i},n}|^{5/2} \leq m_n \left(\frac{2d_4}{d_2}(nh_n^3)^{1/4}\right)^{5/2} \leq d_5 \left(\frac{2d_4}{d_2}\right)^{5/2}(nh_n^3)^{1/8}.$$

This bound when combined with (H) implies that as $n \to \infty$,
$$\sum_{\mathbf{i} \in \mathcal{J}_n} E|X_{\mathbf{i},n}|^{5/2} \to 0,$$

which by the Shergin fact (with $s = 5/2$) yields
$$S_n(\tau) = \sum_{\mathbf{i} \in \mathcal{J}_n} X_{\mathbf{i},n} \xrightarrow{d} Z.$$

Thus, using (2.50) and $\beta_n = P\{X_i \in D_n(\tau + \sqrt{nh_n})\} \to 0$, Lemma 2.2 implies that

(2.55) $$\frac{T_n(\tau)}{\sigma_n(\tau)} \xrightarrow{d} Z.$$

Putting everything together we get from (2.48) that
$$\lim_{\tau \to \infty} \limsup_{n \to \infty} E(R_n(\tau))^2 = 0$$

and from (2.46) that
$$\lim_{\tau \to \infty} \limsup_{n \to \infty} |\sigma_n^2(\tau) - \sigma^2| = 0,$$

which in combination with (2.55) implies that

(2.56) $$a_n(L_n(c) - E\Pi_n(c)) \xrightarrow{d} \sigma Z,$$

where

(2.57) $$\sigma^2 = \int_{-\infty}^{\infty} \int_{I_d} \int_B \Gamma(\theta, s, t) \iota(\theta) \, dt \, d\theta \, ds$$

with $\Gamma(\theta, s, t)$ as defined in (2.28). Since by Lemma 2.3

$$E(a_n(L_n(c) - E\Pi_n(c)))^2 \leq 2 \operatorname{Var}(a_n \Pi_n(c))$$

and

$$\operatorname{Var}(a_n \Pi_n(c)) \to \sigma^2 < \infty,$$

we can conclude that

$$a_n(EL_n(c) - E\Pi_n(c)) \to 0$$

and thus

$$a_n(L_n(c) - EL_n(c)) \xrightarrow{d} \sigma Z.$$

This gives that

(2.58) $$a_n \left( G(\mathbb{C}_n(c) \Delta C(c)) - \int_{E_n} E|I\{f_n(x) \geq c\} - I\{f(x) \geq c\}| \, dG(x) \right)$$
$$\xrightarrow{d} \sigma Z,$$

which is (2.14). In light of (2.58) and keeping mind that

$$EG(\mathbb{C}_n(c) \Delta C(c)) = \int_{\mathbb{R}^d} E|I\{f_n(x) \geq c\} - I\{f(x) \geq c\}|g(x) \, dx,$$

we see that to complete the proof of (2.1) it remains to show that

(2.59) $$a_n E \int_{E_n^c} |I\{f_n(x) \geq c\} - I\{f(x) \geq c\}|g(x) \, dx \to 0.$$

We shall begin by bounding

$$E|I\{f_n(x) \geq c\} - I\{f(x) \geq c\}|, \qquad x \in E_n^c.$$

Applying inequality (2.38) with $a = f_n(x) - f(x)$ and $b = c - f(x)$ we have for $x \in E_n^c$,

$$E|I\{f_n(x) \geq c\} - I\{f(x) \geq c\}|$$
$$\leq EI\{|f_n(x) - f(x)| \geq |c - f(x)|\}$$
$$= P\{|f_n(x) - f(x)| \geq |c - f(x)|\}$$
$$\leq P\{|f_n(x) - Ef_n(x)| \geq |c - f(x)| - |f(x) - Ef_n(x)|\}.$$

By recalling the definition of $E_n$ in (2.17) we obtain

$$E|I\{f_n(x) \geq c\} - I\{f(x) \geq c\}|$$
$$\leq P\bigg\{|f_n(x) - Ef_n(x)| \geq \frac{\varsigma(\log n)^{1/2}}{(nh_n)^{1/2}} - |f(x) - Ef_n(x)|\bigg\}$$
$$\leq P\bigg\{|f_n(x) - Ef_n(x)| \geq \frac{\varsigma(\log n)^{1/2}}{(nh_n)^{1/2}} - A_1 h_n^{2/d}\bigg\}.$$

The last inequality uses the fact that (K.i), (K.ii), (2.7) and (2.8) imply after a change of variables and an application of Taylor's formula for $f(x + h_n^{1/d}v) - f(x)$ that for some constant $A_1 > 0$,

$$\sup_{n \geq 2} h_n^{-2/d} \sup_{x \in \mathbb{R}^d} |Ef_n(x) - f(x)| \leq A_1.$$

Thus for all large enough $n$ uniformly in $x \in E_n^c$,

$$E|I\{f_n(x) \geq c\} - I\{f(x) \geq c\}|$$
$$\leq P\bigg\{|f_n(x) - Ef_n(x)| \geq \frac{\varsigma(\log n)^{1/2}}{(nh_n)^{1/2}}\bigg(1 - \frac{A_1}{\varsigma}\frac{\sqrt{nh_n^{1+2/d}}h_n^{1/d}}{\sqrt{\log n}}\bigg)\bigg\}$$
$$\leq P\bigg\{|f_n(x) - Ef_n(x)| \geq \frac{\varsigma(\log n)^{1/2}}{2(nh_n)^{1/2}}\bigg\} =: p_n(x),$$

where the last inequality uses (H). We shall bound $p_n(x)$ using Bernstein's inequality on the i.i.d. sum

$$f_n(x) - Ef_n(x) = \frac{1}{nh_n}\sum_{i=1}^n \bigg\{K\bigg(\frac{x-X_i}{h_n^{1/d}}\bigg) - EK\bigg(\frac{x-X_i}{h_n^{1/d}}\bigg)\bigg\}.$$

Notice that for each $i = 1, \ldots, n$,

$$\operatorname{Var}\bigg(\frac{1}{nh_n}K\bigg(\frac{x-X_i}{h_n^{1/d}}\bigg)\bigg) \leq \frac{1}{(nh_n)^2}\int_{\mathbb{R}^d} K^2\bigg(\frac{x-y}{h_n^{1/d}}\bigg)f(y)\,dy$$
$$= \frac{1}{n^2 h_n}\int_{\mathbb{R}^d} K^2(u)f(x - h_n^{1/d}u)\,du \leq \frac{\|K\|_2^2 M}{n^2 h_n}$$

and by (K.i),

$$\frac{1}{nh_n}\bigg|K\bigg(\frac{x-X_i}{h_n^{1/d}}\bigg) - EK\bigg(\frac{x-X_i}{h_n^{1/d}}\bigg)\bigg| \leq \frac{2\kappa}{nh_n}.$$

Therefore by Bernstein's inequality [i.e., page 855 of Shorack and Wellner (1986)],

$$p_n(x) \leq 2\exp\bigg(\frac{-\varsigma^2(\log n)/(4nh_n)}{\|K\|_2^2 M/(nh_n) + 2/3\varsigma(\log n)^{1/2}/(2(nh_n)^{1/2})\kappa/(nh_n)}\bigg)$$

$$= 2\exp\left(\frac{-\varsigma^2(\log n)/4}{\|K\|_2^2 M + \kappa\varsigma(\log n)^{1/2}/(3(nh_n)^{1/2})}\right).$$

Hence by (H) and keeping in mind that $\varsigma > \sqrt{2}$ in (2.17), we get for some constant $a > 0$ that for all large enough $n$, uniformly in $x \in E_n^c$, we have the bound

(2.60) $$p_n(x) \leq 2\exp(-\varsigma a \log n).$$

We shall show below that $\lambda(\mathbb{C}_n(c)\Delta C(c)) \leq m < \infty$ for some $0 < m < \infty$. Assuming this to be true, we have the following [similar lines of arguments are used in Rigollet and Vert (2008)]

$$E\int_{E_n^c} |I\{f_n(x) \geq c\} - I\{f(x) \geq c\}|g(x)\,dx$$

$$= E\int_{E_n^c \cap (\mathbb{C}_n(c)\Delta C(c))} |I\{f_n(x) \geq c\} - I\{f(x) \geq c\}|g(x)\,dx$$

$$\leq \sup_{A:\lambda(A)\leq m} E\int_{E_n^c \cap A} |I\{f_n(x) \geq c\} - I\{f(x) \geq c\}g(x)|\,dx$$

$$\leq \sup_{A:\lambda(A)\leq m} \int_{E_n^c \cap A} E|I\{f_n(x) \geq c\} - I\{f(x) \geq c\}|g(x)\,dx$$

$$\leq m\sup_x g(x) \sup_{x \in E_n^c} E|I\{f_n(x) \geq c\} - I\{f(x) \geq c\}|$$

$$\leq m\sup_x g(x) \sup_{x \in E_n^c} p_n(x).$$

With $c_0 = m\sup_x g(x)$ and (2.60) this gives the bound

$$a_n E\int_{E_n^c} |I\{f_n(x) \geq c\} - I\{f(x) \geq c\}|g(x)\,dx$$

$$\leq 2c_0 a_n \exp(-\varsigma a \log n).$$

Clearly by (H), we see that for large enough $\varsigma > 0$

$$a_n \exp(-\varsigma a \log n) \to 0$$

and thus (2.59) follows. It remains to verify that there exists $0 < m < \infty$ with

$$\lambda(\mathbb{C}_n(c)\Delta C(c)) \leq m.$$

Notice that

$$1 \geq \int_{\mathbb{C}_n(c)} f_n(x)\,dx \geq c\lambda(\mathbb{C}_n(c))$$

and
$$1 \geq \int_{C(c)} f(x) \, dx \geq c\lambda(C(c)).$$

Thus
$$\lambda(\mathbb{C}_n(c) \Delta C(c)) \leq 2/c =: m.$$

We see now that the proof of the theorem in the case $k=1$ and $d \geq 2$ is complete.

The proof for the case $k \geq 2$ goes through by an obvious extension of the argument used in the case $k=1$. On account of (B.ii) we can write for large enough $n$

$$\int_{E_n} \Delta_n(x) \, dG(x) = \sum_{j=1}^{n} \int_{E_{j,n}} \Delta_n(x) \, dG(x),$$

where the sets $E_{j,n}$, $j = 1, \ldots, k$, are disjoint and constructed from the $\beta_j$ just as $E_n$ was formed from the boundary set $\beta$ in the proof for the case $k=1$. Therefore by reason of the Poissonization, the summands are independent. Hence the asymptotic normality readily follows as before, where the limiting variance in (2.1) becomes

(2.61) $$\sigma^2 = \sum_{i=1}^{k} c_j^2 \sigma_j^2,$$

where each $\sigma_j^2$ is formed just like (2.57).

2.6. *Proof of the theorem in the case $d=1$.* The case $d=1$ follows along very similar ideas as presented above in the case $d \geq 2$ and is in fact somewhat simpler than the case $d \geq 2$. We therefore skip all the details and only point out that by assumption (B.ii) the boundary set $\beta = \{x \in \mathbb{R} : f(x) = c\}$ consists of $k$ points $z_i$, $i = 1, \ldots, k$. Therefore, the integral over $\theta$ in the definition of $\sigma^2$ in (2.57) has to be replaced by a sum, leading to

(2.62) $$\sigma^2 := \sum_{i=1}^{k} (g^{(1)}(z_i))^2 \int_{-\infty}^{\infty} \int_{-1}^{1} \Gamma(i,s,t)|s|^{2/\gamma_g} \, dt \, ds,$$

where

$\Gamma(i,s,t)$
$$= \mathrm{cov}\left(\left|I\left\{Z_1 \geq -\frac{sf'(z_i)}{\sqrt{c}\|K\|_2}\right\} - I\left\{0 \geq -\frac{sf'(z_i)}{\sqrt{c}\|K\|_2}\right\}\right|,\right.$$
$$\left.\left|I\left\{\rho(t)Z_1 + \sqrt{1-\rho^2(t)}Z_2 \geq -\frac{sf'(z_i)}{\sqrt{c}\|K\|_2}\right\} - I\left\{0 \geq -\frac{sf'(z_i)}{\sqrt{c}\|K\|_2}\right\}\right|\right).$$

We can drop the absolute value sign on $f'(z_i)$ in our definition of $\Gamma(i,s,t)$ for $i = 1, \ldots, k$ and thus $\sigma^2$, since $\rho(t) = \rho(-t)$. $\square$

2.7. *Remarks on the variance and its estimation.* Clearly the variance $\sigma_G^2$ that appears in Theorem 1 does not have a nice closed form and in many situations is not feasible to calculate. Therefore in applications $\sigma_G^2$ will very likely have to be estimated either by simulation or from the data itself. In the latter case, an obvious suggestion is to apply the bootstrap and another is to use the jackknife. A separate investigation is required to verify that these methods work in this setup. [Similarly we may also need to estimate $EG(\mathbb{C}_n(c)\Delta C(c))$.]

Here is a quick and dirty way to estimate $\sigma_G^2$. Let $X_1, \ldots, X_n$ be i.i.d. $f$. Choose a sequence of integers $1 \leq m_n \leq n$, such that $m_n \to \infty$ and $n/m_n \to \infty$. Set $\varsigma_n = [n/m_n]$ and take a random sample of the data $X_1, \ldots, X_n$ of size $m_n \varsigma_n$ and then randomly divide this sample into $\varsigma_n$ disjoint samples of size $m_n$. Let

$$\xi_i = d_G(\mathbb{C}_{m_n}^{(i)}(c)\Delta C(c)) \qquad \text{for } i = 1, \ldots, \varsigma_n,$$

where $\mathbb{C}_{m_n}^{(i)}(c)$ is formed from sample $i$. We propose as our estimator of $\sigma_G^2$, the sample variance of $(\frac{m_n}{h_{m_n}})^{1/4}\xi_i$, $i = 1, \ldots, \varsigma_n$,

$$\left(\frac{m_n}{h_{m_n}}\right)^{1/2} \sum_{i=1}^{\varsigma_n} (\xi_i - \overline{\xi})^2/(\varsigma_n - 1).$$

Under suitable regularity conditions it is routine to show that this is a consistent estimator of $\sigma_G^2$, again the details are beyond the scope of this paper.

*The variance $\sigma_G^2$ under a bivariate normal model.* In order to obtain a better understanding about the expression of the variance we consider it in the following simple bivariate normal example. Assume

$$f(x, y) = \frac{1}{2\pi} \exp\left(-\frac{x^2 + y^2}{2}\right), \qquad (x, y) \in \mathbb{R}^2.$$

A special case of Theorem 1 says that whenever

(2.63) $$nh_n^2 \to 0 \quad \text{and} \quad nh_n/\log n \to \infty$$

(here $\gamma = 0$), then

(2.64) $$\left(\frac{n}{h_n}\right)^{1/4} \{\lambda(\mathbb{C}_n(c)\Delta C(c)) - E\lambda(\mathbb{C}_n(c)\Delta C(c))\} \xrightarrow{d} \sigma_\lambda Z.$$

We shall calculate $\sigma_\lambda$ in this case. We get that

$$f'(x, y) = -(x, y)f(x, y).$$

Notice for any $0 < c < \frac{1}{2\pi}$,

$$\beta = \{(x, y) : x^2 + y^2 = -2\log(c2\pi)\}.$$

Setting
$$r(c) = \sqrt{-2\log(c2\pi)},$$
we see that $\beta$ is the circle with center 0 and radius $r(c)$. Choosing the obvious differmorphism,
$$y(\theta) = (r(c)\cos\theta, r(c)\sin\theta) \quad \text{for } \theta \in [0, 2\pi],$$
we get that for $\theta \in [0, 2\pi]$,
$$u(\theta) = (-\cos\theta, -\sin\theta), \qquad y'(\theta) = (r(c)\sin\theta, -r(c)\cos\theta)$$
and
$$\iota(\theta) = \left|\det\begin{vmatrix} r(c)\sin\theta & -r(c)\cos\theta \\ -\cos\theta & -\sin\theta \end{vmatrix}\right| = r(c).$$

Here $g = 1$ and we are assuming $\gamma = 0$. We get that
$$c(s, \theta, \gamma t) = c(s, \theta, 0) = -\frac{s|f'(y(\theta))|}{\sqrt{c}\|K\|_2} = -\frac{sr(c)\sqrt{c}}{\|K\|_2}.$$

Thus
$$\Gamma(\theta, s, t)$$
$$= \text{cov}\left(\left|I\left\{Z_1 \geq -\frac{sr(c)\sqrt{c}}{\|K\|_2}\right\} - I\left\{0 \geq -\frac{sr(c)\sqrt{c}}{\|K\|_2}\right\}\right|,\right.$$
$$\left.\left|I\left\{\rho(t)Z_1 + \sqrt{1-\rho^2(t)}Z_2 \geq -\frac{sr(c)\sqrt{c}}{\|K\|_2}\right\} - I\left\{0 \geq -\frac{sr(c)\sqrt{c}}{\|K\|_2}\right\}\right|\right).$$

This gives
$$\sigma_\lambda^2 = r(c)\int_{-\infty}^{\infty}\int_0^{2\pi}\int_B \Gamma(\theta, s, t)\, dt\, d\theta\, ds.$$

Set
$$\Upsilon(\theta, u, t) = \text{cov}(|I\{Z_1 \geq -u\} - I\{0 \geq -u\}|,$$
$$|I\{\rho(t)Z_1 + \sqrt{1-\rho^2(t)}Z_2 \geq -u\} - I\{0 \geq -u\}|).$$

We see then by the change of variables $u = \frac{sr(c)\sqrt{c}}{\|K\|_2}$ that
$$\sigma_\lambda^2 = \frac{\|K\|_2}{\sqrt{c}}\int_{-\infty}^{\infty}\int_0^{2\pi}\int_B \Upsilon(\theta, u, t)\, dt\, d\theta\, du.$$

For comparison, Theorem 2.1 of Cadre (2006) says that if
(2.65) $\qquad nh_n/(\log n)^{16} \to \infty \quad \text{and} \quad nh_n^3(\log n)^2 \to 0,$

then

$$(2.66) \quad \sqrt{nh_n}\lambda(\mathbb{C}_n(c)\Delta C(c)) \xrightarrow{P} \|K\|_2\sqrt{\frac{2c}{\pi}}\int_\beta \frac{d\mathcal{H}}{\|\nabla f\|} = 2\|K\|_2\sqrt{\frac{2\pi}{c}}.$$

The measure $d\mathcal{H}$ denotes the Hausdorff measure on $\beta$. In this case $\mathcal{H}(\beta)$ is the circumference of $\beta$.

Observe that since $\sqrt{nh_n}(\frac{h_n}{n})^{1/4} = (nh_n^2)^{1/4}h_n^{1/4} \to 0$, (2.64) and (2.66) imply that whenever (2.63) and (2.65) hold, we get

$$\sqrt{nh_n}E\lambda(\mathbb{C}_n(c)\Delta C(c)) \to 2\|K\|_2\sqrt{\frac{2\pi}{c}}.$$

[Notice that the choice $h_n = 1/(\sqrt{n\log n})$ satisfies both (2.63) and (2.65).]

**Acknowledgments.** The authors would like to thank Gerard Biau for pointing out to them the *tubular neighborhood theorem*, and also the Associate Editor and the referee for useful suggestions and a careful reading of the manuscript.

## REFERENCES

BAÍLLO, A., CUEVAS, A. and JUSTEL, A. (2000). Set estimation and nonparametric detection. *Canad. J. Statist.* **28** 765–782. MR1821433

BAÍLLO, A., CUESTAS-ALBERTOS, J. A. and CUEVAS, A. (2001). Convergence rates in nonparametric estimation of level sets. *Statist. Probab. Lett.* **53** 27–35. MR1843338

BAÍLLO, A. (2003). Total error in a plug-in estimator of level sets. *Statist. Probab. Lett.* **65** 411–417. MR2039885

BAÍLLO, A. and CUEVAS, A. (2006). Image estimators based on marked bins. *Statistics* **40** 277–288. MR2263620

BEIRLANT, J. and MASON, D. M. (1995). On the asymptotic normality of $L_p$-norms of empirical functionals. *Math. Methods Statist.* **4** 1–19. MR1324687

BIAU, G., CADRE, B. and PELLETIER, B. (2008). Exact rates in density support estimation. *J. Multivariate Anal.* **99** 2185–2207.

BREDON, G. E. (1993). *Topology and Geometry. Graduate Texts in Mathematics* **139**. Springer, New York. MR1224675

CADRE, B. (2006). Kernel estimation of density level sets. *J. Multivariate Anal.* **97** 999–1023. MR2256570

CAVALIER, L. (1997). Nonparametric estimation of regression level sets. *Statistics* **29** 131–160. MR1484386

CUEVAS, A., FEBRERO, M. and FRAIMAN, R. (2000). Estimating the number of clusters. *Canad. J. Statist.* **28** 367–382.

CUEVAS, A., GONZÁLEZ-MANTEIGA, W. and RODRÍGUEZ-CASAL, A. (2006). Plug-in estimation of general level sets. *Aust. N. Z. J. Statist.* **48** 7–19. MR2234775

DESFORGES, M. J., JACOB, P. J. and COOPER, J. E. (1998). Application of probability density estimation to the detection of abnormal conditions in engineering. *Proceedings of the Institute of Mechanical Engineering* **212** 687–703.

EINMAHL, U. and MASON, D. M. (2005). Uniform in bandwidth consistency of kernel-type function estimators. *Ann. Statist.* **33** 1380–1403. MR2195639

Fan, W., Miller, M., Stolfo, S. J., Lee, W. and Chan, P. K. (2001). Using artificial anomalies to detect unknown and known network intrusions. In *IEEE International Conference on Data Mining (ICDM'01)* 123–130. IEEE Computer Society.

Feller, W. (1966). *An Introduction to Probability Theory and Its Applications, Vol. II.* Wiley, New York. MR0210154

Gardner, A. B., Krieger, A. M., Vachtsevanos, G. and Litt, B. (2006). One-class novelty detection for seizure analysis from intracranial EEG. *J. Mach. Learn. Res.* **7** 1025–1044. MR2274396

Gayraud, G. and Rousseau, J. (2005). Rates of convergence for a Bayesian level set estimation. *Scand. J. Statist.* **32** 639–660. MR2232347

Gerig, G., Jomier, M., Chakos, M. (2001). VALMET: A new validation tool for assessing and improving 3D object segmentation. In *Medical Image Computing and Computer Assisted Intervention MICCAI* **2208** (W. Niessen and M. Viergever, eds.) 516–523. Springer, New York.

Giné, E., Mason, D. M. and Zaitsev, A. Y. (2003). The $L_1$-norm density estimator process. *Ann. Probab.* **31** 719–768. MR1964947

Goldenshluger, A. and Zeevi, A. (2004). The Hough transform estimator. *Ann. Statist.* **32** 1908–1932. MR2102497

Hall, P. and Kang, K.-H. (2005). Bandwidth choice for nonparametric classification. *Ann. Statist.* **33** 284–306. MR2157804

Hartigan, J. A. (1975). *Clustering Algorithms*. Wiley, New York. MR0405726

Hartigan, J. A. (1987). Estimation of a convex density contour in two dimensions. *J. Amer. Statist. Assoc.* **82** 267–270. MR883354

Horváth, L. (1991). On $L_p$-norms of multivariate density estimators. *Ann. Statist.* **19** 1933–1949. MR1135157

Huo, X. and Lu, J.-C. (2004). A network flow approach in finding maximum likelihood estimate of high concentration regions. *Comput. Statist. Data Anal.* **46** 33–56. MR2056823

Jang, W. (2006). Nonparametric density estimation and clustering in astronomical sky surveys. *Comput. Statist. Data Anal.* **50** 760–774. MR2207006

Johnson, W. B., Schechtman, G. and Zinn, J. (1985). Best constants in moment inequalities for linear combinations of independent and exchangeable random variables. *Ann. Probab.* **13** 234–253. MR770640

King, S. P., King, D. M., Anuzis, P., Astley, L., Tarassenko, K., Hayton, P. and Utete, S. (2002). The use of novelty detection techniques for monitoring high-integrity plant. In *Proceedings of the 2002 International Conference on Control Applications* **1** 221–226. Cancun, Mexico.

Klemelä, J. (2004). Visualization of multivariate density estimates with level set trees. *J. Comput. Graph. Statist.* **13** 599–620. MR2087717

Klemelä, J. (2006a). Visualization of multivariate density estimates with shape trees. *J. Comput. Graph. Statist.* **15** 372–397. MR2256150

Klemelä, J. (2008). Visualization of scales of multivariate density estimates. Unpublished manuscript.

Lang, S. (1997). *Undergraduate Analysis*, 2nd ed. Springer, New York. MR1476913

Mason, D. and Polonik, W. (2008). Asymptotic normality of plug-in level set estimates (expanded version). Unpublished manuscript.

Markou, M. and Singh, S. (2003). Novelty detection: A review—part 1: Statistical approaches. *Signal Processing* **83** 2481–2497.

Molchanov, I. S. (1998). A limit theorem for solutions of inequalities. *Scand. J. Statist.* **25** 235–242. MR1614288

Nairac, A., Townsend, N., Carr, R., King, S., Cowley, L. and Tarassenko, L. (1997). A system for the analysis of jet engine vibration data. *Integrated Comput. Aided Eng.* **6** 53–65.

Polonik, W. (1995). Measuring mass concentrations and estimating density contour clusters—an excess mass approach. *Ann. Statist.* **23** 855–881. MR1345204

Prastawa, M., Bullitt, E., Ho, S. and Gerig, G. (2003). Robust estimation for brain tumor segmentation. In *MICCAI Proceedings LNCS 2879* 530–537. Springer, Berlin.

Rigollet, P. and Vert, R. (2008). Fast rates for plug-in estimators of density level sets. Available at arxiv.org/pdf/math/0611473.

Roederer, M. and Hardy, R. R. (2001). Frequency difference gating: A multivariate method for identifying subsets that differ between samples. *Cytometry* **45** 56–64.

Rosenblatt, M. (1975). A quadratic measure of deviation of two-dimensional density estimates and a test of independence. *Ann. Statist.* **3** 1–14. MR0428579

Scott, C. D. and Davenport, M. (2006). Regression level set estimation via cost-sensitive classification. *IEEE Trans. Inf. Theory.* **55** 2752–2757.

Scott, C. D. and Nowak, R. D. (2006). Learning minimum volume sets. *J. Mach. Learn. Res.* **7** 665–704. MR2274383

Shergin, V. V. (1990). The central limit theorem for finitely-dependent random variables. In *Probability Theory and Mathematical Statistics, Vol. II (Vilnius, 1989)* 424–431. "Mokslas," Vilnius. MR1153895

Shorack, G. R. and Wellner, J. A. (1986). *Empirical Processes with Applications to Statistics*. Wiley, New York. MR838963

Steinwart, I., Hush, D. and Scovel, C. (2004). Density level detection is classification. Technical report, Los Alamos national laboratory.

Steinwart, I., Hush, D. and Scovel, C. (2005). A classification framework for anomaly detection. *J. Mach. Learn. Res.* **6** 211–232 (electronic). MR2249820

Stuetzle, W. (2003). Estimating the cluster type of a density by analyzing the minimal spanning tree of a sample. *J. Classification* **20** 25–47. MR1983120

Theiler, J. and Cai, D. M. (2003). Resampling approach for anomaly detection in multispectral images. In *Proceedings of the SPIE 5093* 230–240.

Tsybakov, A. B. (1997). On nonparametric estimation of density level sets. *Ann. Statist.* **25** 948–969. MR1447735

Tsybakov, A. B. (2004). Optimal aggregation of classifiers in statistical learning. *Ann. Statist.* **32** 135–166. MR2051002

Vert, R. and Vert, J.-P. (2006). Consistency and convergence rates of one-class SVMs and related algorithms. *J. Mach. Learn. Res.* **7** 817–854. MR2274388

Walther, G. (1997). Granulometric smoothing. *Ann. Statist.* **25** 2273–2299. MR1604445

Wand, M. (2005). Statistical methods for flow cytometric data. Presentation. Available at http://www.uow.edu.au/~mwand/talks.html.

Willett, R. M. and Nowak, R. D. (2005). Level set estimation in medical imaging. In *Proceedings of the IEEE Statistical Signal Processing Workshop*. Bordeaux, France.

Willett, R. M. and Nowak, R. D. (2006). Minimax optimal level set estimation. Submitted to *IEEE Trans. Image Proc.* Available at http://www.ee.duke.edu/~willett/.

Yeung, D. W. and Chow, C. (2002). Parzen window network intrusion detectors. In *Proceedings of the 16th International Conference on Pattern Recognition* **4** 385–388. Quebec, Canada.

Statistics Program  
University of Delaware  
206 Townsend Hall  
Newark, Delaware 19717  
USA  
E-mail: davidm@udel.edu

Departement of Statistics  
University of California, Davis  
One Shields Avenue  
Davis, California 95616–8705  
USA  
E-mail: wpolonik@ucdavis.edu